*Research Article - Preprint*

## A Trust-region Funnel Algorithm for Grey-Box Optimisation

Gul Hameed [1], Tao Chen [1], Antonio del Rio Chanona [2], Lorenz T Biegler [3], Michael Short [1*]

[1] School of Chemistry and Chemical Engineering, University of Surrey, Guildford GU2 7XH, UK
[2] Department of Chemical Engineering, Imperial College London, United Kingdom
[3] Department of Chemical Engineering, Carnegie Mellon University, Pittsburgh, PA 15213, USA

*Corresponding Author: m.short@surrey.ac.uk

**Abstract**

Grey-box optimisation, where some parts of an optimisation problem are represented by explicit algebraic (glass-box) models while others are treated as black-box models lacking analytic derivatives, remains a challenge in process systems engineering. Trust-region (TR) methods provide a robust framework for grey-box problems by combining accurate glass-box derivatives with local reduced models (RMs) for black-box components. However, existing TR approaches often involve complex multi-layered formulations requiring extensive parameter tuning, or lack open-source implementations. Motivated by the recent advances in funnel-based convergence theory for nonlinear optimisation and the TR filter method, we propose a novel TR funnel algorithm for grey-box optimisation that replaces the filter acceptance criterion with a generalisable uni-dimensional funnel, maintaining a monotonically non-increasing upper bound (a funnel) on approximation error of the local black-box RMs. A global convergence proof to a first-order critical point is established. The algorithm, implemented in an open-source Pyomo framework, supports multiple RM forms and globalisation strategies (filter or funnel). Benchmark tests on seven numerical and engineering problems show that the TR funnel algorithm achieves comparable and often improved performance relative to the classical TR filter method. The TR funnel method thus provides a simpler, and extensible alternative for large-scale grey-box optimisation.

# 1. Introduction

Process optimisation has advanced significantly in recent decades, enabling simultaneous simulation and optimisation of increasingly complex process flowsheets. Early sequential modular simulators solved unit operations in order of appearance but often suffered from noisy derivative estimates, causing gradient-based solvers to fail or deviate from optimal solutions, particularly in large-scale systems with many degrees of freedom [1]. Modern equation-oriented (EO) frameworks address these limitations by treating the entire process as a coupled system of nonlinear equations and computing exact derivatives via automatic differentiation [1]. This shift has enabled powerful nonlinear programming (NLP) solvers based on active-set (CONOPT [2], MINOS [3], SNOPT [4]) and primal-dual interior point methods (KNITRO [5], LOQO [6], IPOPT [7]) to compute locally optimal solutions.

Many industrial optimisation problems, however, lack complete EO representation. Certain process units rely on computationally expensive or proprietary models (e.g., computational fluid dynamics simulations or external simulators such as Aspen) that can only be evaluated as black-box functions with no or limited derivative information. Constructing global surrogates that are sufficiently accurate for such units can be impractical due to high computational cost of sampling [9]. A more efficient strategy is to build local surrogates, also known as reduced models (RMs), around the current iterate, approximating the black-box response without extensive sampling. This results in a grey-box optimisation problem, where some model components are explicit (glass-box) equations with known derivatives, while others are black-box functions lacking analytic gradients (i.e., without mathematical formulations). An effective grey-box solver exploits the cheap, accurate derivatives of the glass-box model while simultaneously reducing both objective function and black-box approximation error (measured as an infeasibility measure).

RM-based trust-region (TR) frameworks are well-suited to this challenge. At each iteration, they solve a local TR subproblem (TRSP) that combines the glass-box model with local RMs for black-box components, within a bounded region around the current iterate. The region size is dynamically adjusted based on RM prediction accuracy. The TR filter method [1][9][10] has emerged as a robust approach in process systems engineering, balancing the infeasibility measure and optimality in grey-box problems. The TR filter method evaluates a single trial step by solving a TRSP in each iteration. It maintains a list of pairs (objective, infeasibility measure) from past iterates (known as the filter set) and accepts a trial step if it achieves a sufficient decrease in either the objective or the infeasibility measure without worsening the other. Constraint violations are handled implicitly by the NLP solver used to solve the TRSP (e.g., IPOPT). In practice, the filter approach is effective but requires the manual adjustment of the tuning parameters and may stall when RM errors dominate the optimisation process.

An alternative is the derivative-free trust-funnel (DEFT-FUNNEL) algorithm [11], which enforces a monotonically decreasing bound ("funnel") on constraint violation (capturing both glass-box and black-box effects). Unlike the TR Filter method which is a single-run algorithm converging to a local optimum, the DEFT-FUNNEL method attempts global optimality by combining a global multi-start search (Multi-Level Single Linkage, MLSL) and a trust-funnel sequential quadratic optimisation (SQO) local search. Each SQO iteration computes a composite step consisting of a normal step (reducing constraint violation) and a tangential step (improving the objective without worsening the constraint violation). Black-box functions are approximated via polynomial interpolation RMs, and a self-correcting geometry approach [12] manages the interpolation set. While promising, the DEFT-FUNNEL is algorithmically complex due to the integration of MLSL, the funnel, SQO composite steps and geometry management, and its MATLAB implementation limits accessibility compared to the open-source Pyomo-based TR Filter method.

ARGONAUT is another derivative-free/grey-box optimiser [13], in which the black-box components are replaced by RMs. A global optimisation algorithm explores the design space to provide a lower bound while a local optimisation with different starting points establishes upper bounds. After solving, the RMs are updated using function evaluations collected from the global and local minimum found within a clustering procedure that selects new sample points. This process of solving, updating RMs, and re-solving is repeated until convergence. Another grey-box optimisation approach, two-phase TR algorithm [14], employs radial basis function (RBF) RMs to approximate the black-box functions. It first runs a feasibility phase to locate a feasible solution, followed by an optimisation phase that searches for a global minimum starting from the feasible point. While ARGONAUT and two-phase algorithms have been successfully applied to carbon capture optimisation [14][15], their adoption in process optimisation remains limited due to a lack of open-source availability.

Despite the many advances, challenges remain in the development of efficient grey-box and derivative-free optimisation (DFO) methods [16]. Few algorithms provide rigorous convergence guarantees within finite iterations. Sampling requirements increase rapidly with black-box dimensionality, limiting scalability. Problems with many grey-/black-box constraints are particularly challenging, as most frameworks cannot efficiently handle large grey-box formulations without user intervention (e.g., manually adjusting the tuning parameters). Finally, there is a need for algorithms that exploit the complementary strengths of glass-box and black-box models, maximising communication between the two to achieve robust convergence.

Among available grey-box optimisation approaches, the classical TR filter method and its extensions remain the most widely applied, owing to their algorithmic simplicity and demonstrated success in diverse applications. Use cases include the optimisation of the Williams-Otto process, ammonia synthesis, air-fired power-plant [1], $CO_2$ capture process flowsheets [9], pressure swing adsorption (PSA) process [17][18], heat and mass exchanger networks [19][20], chlorobenzene process [21], refinery process [22], carbon capture plant design [23] and produced water networks [24]. Recent advancements focused on enhancing reduced-order modelling beyond polynomial interpolation (using first-order Taylor expansions [21][17] and Gaussian Process regression [25]) and integrating local Hessian information [26] to accelerate convergence. Nevertheless, a recent review on RM-based TR filter methods [27] highlights challenges: high numbers of expensive black-box evaluations and sensitivity to algorithmic tuning parameters.

Recently, the Unified Nonlinear Optimisation (UNO) solver [28] reinterpreted the funnel method through the lens of a filter framework, providing a rigorous global convergence proof (initially derived for the filter method) for the funnel approach. Motivated by this development in NLP and the classical TR filter method, we propose a new TR funnel method for grey-box optimisation and derive its global convergence proof to local optimality. Our approach replaces the filter acceptance criterion with a dynamically updated funnel constraint on the infeasibility measure, reducing sensitivity to the tuning parameters, retaining the robust local convergence guarantees of the classical TR filter framework, and maintaining comparable performance. Unlike the filter method, which maintains and updates a list of non-dominated (objective, infeasibility) pairs, the funnel only updates a single scalar parameter "the funnel width" (i.e., the bound on the infeasibility measure) at each iteration. We implement the TR funnel algorithm in an open-source Pyomo framework, allowing users to select among multiple RM forms (linear, quadratic, Taylor expansions, Gaussian Process regression) [26] and globalisation strategies (filter or funnel). Benchmark tests on challenging numerical and engineering problems demonstrate that the TR funnel method matches the performance of the classical TR filter algorithm while offering a simpler and more extensible framework for grey-box optimisation.

The remainder of this paper is organised as follows. Section 2 presents the formulation of the proposed TR funnel algorithm. Section 3 (supported by Appendix A) establish the convergence proof to a first-order critical point. Section 4 describes the algorithmic implementation in open-source Pyomo. Section 5 reports numerical and engineering case tests results. Finally, Section 6 summarises the key findings and outlines the future research directions.

## 2. Trust-Region Funnel Framework

The TR framework has long been used to solve NLPs by constructing local approximations of the objective and constraint functions and iteratively refining them within dynamically adjusted TRs [29][30]. In grey-box NLP optimisation, the TR approach additionally controls approximation errors introduced by RMs built for expensive black-box components [18][31].

We consider a general grey-box optimisation problem of the form:

$$\min_{w,z} f\big(z, w, t(w)\big) \ s.t. \ \ h\big(z, w, t(w)\big) = 0 \ and \ g(z, w, t(w)) \leq 0, \qquad (1)$$

where $z \in \mathbb{R}^n$ and $w \in \mathbb{R}^m$ are vector-valued decision variables, and the black-box mapping $t(w) : \mathbb{R}^m \to \mathbb{R}^p$ provides a vector of additional process outputs that are expensive to evaluate or lack analytic derivatives. The functions $f$, $h$, $g$ and $t$ are assumed twice continuously differentiable.

To decouple the glass-box and black-box components, an auxiliary variable vector $y \in \mathbb{R}^p$ is introduced, resulting in an equivalent reformulation:

$$\min_x f(x) \quad s.t. \quad h(x) = 0, g(x) \leq 0 \; and \; y = t(w). \tag{2}$$

Vector $x^T = [w^T, y^T, z^T]$ denotes decision variables in all formulations. This formulation separates explicit model equations from implicit black-box evaluations, allowing local reduced models to be constructed only for $t(w)$.

In DFO, the TR typically serves two purposes: limiting the step size and defining the region where RMs remain accurate. In grey-box settings, however, these two roles can conflict, as applying the same region to both black- and glass-box variables may unnecessarily restrict progress. To resolve this, a separate sampling region has previously been introduced [9][26][32], used solely for RM construction, while the TR continues to govern step acceptance and globalisation.

At each iteration $k$, a sampling region $S(w^{(k)}, \sigma^{(k)})$ is defined around the current black-box inputs $w^{(k)}$, with radius $\sigma^{(k)} \leq \Delta^{(k)}$, ensuring it is contained within the TR radius $\Delta^{(k)}$. RMs $r^{(k)}(w)$ are constructed using black-box evaluations within the sampling region, satisfying the *κ-fully linear* property:

$$\left\|\nabla r^{(k)}(w) - \nabla t(w)\right\| \leq \kappa_g \sigma^{(k)} \; and \; \left\|r^{(k)}(w) - t(w)\right\| \leq \kappa_h {\sigma^{(k)}}^2 \tag{3}$$

for all $w$: $\left\|w - w^{(k)}\right\| \leq \Delta^{(k)}$, where $\kappa_g, \kappa_h > 0$. The existence of $\nabla t(w)$ is assumed, even if unavailable. The *κ-fully linear* property ensures that the RMs accurately approximate the true black-box functions $t(w)$ and their gradients $\nabla t(w)$ within the sampling region [16]. Throughout the paper, superscript $k$ indicates iteration number and $\|.\|$ is Euclidean norm unless specified otherwise.

Unlike traditional TR-based methods for NLPs where $\Delta^{(k)} \rightarrow 0$ ensures convergence, in the TR approach used to solve local nonlinear subproblem in this work, only $\sigma^{(k)} \rightarrow 0$ necessarily contracts iteratively (see global convergence proof in Appendix A). This allows larger TR steps for decoupled glass-box variables even near optimality. If condition (3) is not met, RM is re-constructed before proceeding.

The TRSP solved at iteration $k$ is:

$$\min_x f(x) \quad s.t. \; h(x) = 0, g(x) \leq 0, y = r^{(k)}(w) \; and \; \left\|x - x^{(k)}\right\| \leq \Delta^{(k)}. \tag{4}$$

This formulation guarantees feasibility for the explicit constraints, consistency between $y$ and the RM prediction, and restricts the search to a region where RM is reliable. The funnel framework then governs how this subproblem is solved and iterates are advanced. The following sections explain criticality (Section 2.1), compatibility (Section 2.2), TRSP solution (Section 2.3), funnel and TR update (Section 2.4), and feasibility restoration (Section 2.5) mechanisms to ensure accuracy and convergence while minimising expensive black-box evaluations.

## 2.1 Criticality Check

At each iteration $k$, the criticality measure $\chi^{(k)}$ is computed to assess proximity to a first-order Karush–Kuhn–Tucker (KKT) point of the grey-box optimisation problem. The measure is obtained by solving the following linearised problem at the current iterate (using a polyhedral norm such as $l_\infty$ or $l_1$ in the step constraint):

$$\chi^{(k)} = \left|\min_v \nabla f\left(x^{(k)}\right)^T v\right| \tag{5}$$

$$s.t. \nabla h\left(x^{(k)}\right)^T v = 0, g\left(x^{(k)}\right) + \nabla g\left(x^{(k)}\right)^T v \leq 0, v_y - \nabla r^{(k)}\left(w^{(k)}\right)^T v_w = 0, \|v\| \leq 1,$$

where $v^T = \left[v_w^T, v_y^T, v_z^T\right]$ is a normalised trial step.

A small $\chi^{(k)}$ indicates that no significant feasible descent direction exists within the current sampling region, suggesting proximity to stationarity. To ensure that RMs remain accurate, the criticality test is defined relative to the sampling radius:

$$\chi^{(k)} < \xi \sigma^{(k)}, \qquad \xi > 0, \tag{6}$$

and, when satisfied, the sampling radius is reduced:

$$\sigma^{(k)} = \max\left(\min\left(\sigma^{(k-1)}, \chi^{(k)}/\xi\right), \Delta_{min}\right). \tag{7}$$

Progressive contraction $\sigma^{(k)} \to 0$, a small $\chi^{(k)} \to 0$ and the $\kappa - fully\ linear$ property guarantee that limit points generated by the algorithm satisfy first-order KKT conditions of the original grey-box problem [16]. Here, the $\kappa$-*fully linear* property (3) ensure that RM errors are directly proportional to the sampling radius $\sigma^{(k)}$ [16]. As $\sigma^{(k)}$ decreases iteratively, the RM error reduces progressively, resulting in increasingly accurate black-box approximations. This error reduction is essential for reliable descent directions and convergence to KKT conditions.

### 2.2 Compatibility Check

The compatibility check ensures that the representative TRSP remains feasible after introducing the TR constraint and replacing the true black-box response $t(w)$ with its RM approximation $r^{(k)}(w)$ in the grey-box problem formulation. At iteration $k$, the feasibility of the TRSP depends on whether the RM output is sufficiently consistent with the true process constraints within the local TR. This is assessed by solving:

$$\min_{w,z} \left\|y - r^{(k)}(w)\right\| \quad s.t. \;\; h(x) = 0, g(x) \leq 0 \; and \; \left\|x - x^{(k)}\right\| \leq \kappa_\Delta \Delta^{(k)} min\left[1, \kappa_\mu {\Delta^{(k)}}^\mu\right], \tag{8}$$

where $y$ is the fixed RM output, and $\kappa_\Delta \in (0,1)$, $\kappa_\mu > 0$ and $\mu \in (0,1)$ are constants controlling the feasibility search region. Problem (8) is always feasible at $x = x^{(k)}$, where the process constraints are satisfied by construction.

Define $\alpha^{(k)} = \left\|y - r^{(k)}(w)\right\|$. If $\alpha^{(k)} \leq \epsilon_{comp}$, then TRSP is considered compatible, meaning that RM and process constraints yield a feasible local model for step computation. Otherwise (i.e., $\alpha^{(k)} > \epsilon_{comp}$), a restoration phase is initiated to refine the RM and adjust the TR radius for the next iteration.

Upon passing the compatibility check, a feasible step

$$d^{(k)} = x - x^{(k)}, \tag{9}$$

is computed, where $d^{(k)} = \left[d_w^{(k)}, d_y^{(k)}, d_z^{(k)}\right]^T$.

### 2.3 Trust-Region Subproblem

The TRSP is initialised at $x^{(k)} + d^{(k)}$ and solved to obtain a candidate solution $x_s^{(k)}$. The resulting trial step is defined as:

$$s^{(k)} = x_s^{(k)} - x^{(k)}. \tag{10}$$

Analogous to tangential/feasibility and optimality steps in traditional composite-step NLP solvers, $d^{(k)}$ drives the grey-box optimisation to feasibility, while $s^{(k)}$ improves the objective value.

## 2.4 Funnel Mechanism and Trust-Region Update

In TR funnel methods, the infeasibility measure

$$\theta(x) = \|r(w) - t(w)\| \tag{11}$$

evaluates the local consistency between RM predictions and external black-box response.

Several NLP-based strategies have been proposed to manage infeasibility in grey-box optimisation [18]. Penalty methods add a weighted infeasibility term to the objective function, requiring careful tuning of the penalty parameter. Filter methods, on the other hand, rely on heuristic filter parameters, and treat the objective function and infeasibility measure as competing objectives. The funnel method avoids both penalty and filter parameters by maintaining a progressively tightening "funnel-shaped" region that controls allowable infeasibility measure throughout the iterations [28][33].

At iteration $k$, a necessary funnel condition

$$\theta\left(x^{(k)} + s^{(k)}\right) \leq \phi^{(k)} \tag{12}$$

is checked for the acceptance of the trial step, where $\phi^{(k)} > 0$ denotes the current funnel width. The funnel width is initialised as

$$\phi^{(0)} := max\left\{\phi_{min}, \kappa_{\phi}. \theta(x^{(0)})\right\}, \tag{13}$$

with constants $\phi_{min} > 0$ (i.e., hard floor for the funnel width) and $\kappa_{\phi} > 0$ ensuring that the starting point is feasible relative to the initial tolerance. For all $k$, $\phi^{(k+1)} \leq \phi^{(k)}$. This property combined with a finite termination of the feasibility restoration phase guarantees convergence toward a feasible solution [28].

**Iteration Types and Switching Condition**

A switching condition

$$f\left(x^{(k)}\right) - f\left(x^{(k)} + s^{(k)}\right) \geq \delta(\theta\left(x^{(k)}\right))^{\gamma_s}, \qquad \delta \in (0,1), \qquad \gamma_s > 1/(1+\mu) \tag{14}$$

is checked for the step satisfying funnel condition (12), distinguishing between candidate $f - type$ (that improve optimality) and candidate $\theta - type$ steps (that reduce the infeasibility measure).

- $\boldsymbol{f - type}$ **step:** If the switching condition (14) holds, the step is accepted ($x^{(k+1)} = x^{(k)} + s^{(k)}$) if it satisfies a simplified Armijo-type sufficient decrease condition (15).

$$f\left(x^{(k)}\right) - f\left(x^{(k)} + s^{(k)}\right) \geq \eta\Delta^{(k)}, \qquad \eta \in (0,1) \tag{15}$$

The TR radius is enlarged using $\Delta^{(k+1)} := \max\,[\gamma_e\left\|s^{(k)}\right\|, \Delta^{(k)}]$, with $\gamma_e > 1$. Higher $\eta$ slows down convergence to optimality but still ensures feasibility via $\theta - type$ steps.

- **$\theta - type$ step:** If the switching condition (14) is violated, infeasibility is reduced. The step is accepted ($x^{(k+1)} = x^{(k)} + s^{(k)}$) if the funnel sufficient decrease condition (16) is satisfied.

$$\theta\left(x^{(k)} + s^{(k)}\right) \leq \tau\phi^{(k)}, \qquad \tau \in (0,1) \tag{16}$$

  The funnel width is reduced:

$$\phi^{(k+1)} = \left(1 - \kappa_f\right)\theta\left(x^{(k)} + s^{(k)}\right) + \kappa_f\phi^{(k)}, \qquad \kappa_f \in (0, 1). \tag{17}$$

  If $\kappa_f$ (slope-shrink factor for the funnel) is large, the slow decrease in funnel width reduces the progress to feasibility, leading to degraded performance. The TR size is updated using:

$$\Delta^{(k+1)} = \begin{cases} \gamma_c\left\|s^{(k)}\right\| & if\ \rho^{(k)} < \eta_1, \\ \Delta^{(k)} & if\ \eta_1 \leq \rho^{(k)} < \eta_2, \\ \max\left[\gamma_e\left\|s^{(k)}\right\|, \Delta^{(k)}\right] & if\ \rho^{(k)} \geq \eta_2, \end{cases} \tag{18}$$

  where $0 < \gamma_c < 1 < \gamma_e$ and $0 < \eta_1 \leq \eta_2 < 1$. Reduction ratio $\rho^{(k)}$ is calculated as:

$$\rho^{(k)} = \frac{\theta(x^{(k)}) - \theta(x^{(k)} + s^{(k)}) + \epsilon_\theta}{\max\left(\left(\theta_r^{(k)}(x^{(k)}) - \theta_r^{(k)}(x^{(k)} + s^{(k)})\right), \epsilon_\theta\right)} = \frac{\theta(x^{(k)}) - \theta(x^{(k)} + s^{(k)}) + \epsilon_\theta}{\max\left(\left\|t(w^{(k)}) - r^{(k)}(w^{(k)})\right\|, \epsilon_\theta\right)}, \tag{19}$$

  where $\epsilon_\theta > 0$ is a small tolerance to ensure numerical stability near $\theta = 0$ region.
- **Rejected step:** If the step is neither $f - type$ nor $\theta - type$, it is rejected: $x^{(k+1)} = x^{(k)}$, and the TR radius is reduced: $\Delta^{(k+1)} = \gamma_c\left\|s^{(k)}\right\|$.

The funnel mechanism is illustrated in Figure 1. In Figure 1(a), the solid vertical blue lines represent the funnel width, the dotted line indicates funnel envelope, and the black dot represents current iterate. Green and red dots correspond to possible accepted and rejected trial points, respectively, based on the funnel conditions. The $\theta - type$ is presented in Figure 1(b) where the switching condition is violated but the funnel sufficient decrease condition is satisfied. The trial point is accepted, and the funnel width is reduced according to (17), ensuring that the feasibility is driven to zero. The previous funnel width is marked as light blue line. In Figure 1(c), an $f - type$ iteration is shown where both the switching condition and the Armijo-type sufficient decrease condition hold. The trial iterate is accepted, but the funnel width remains unchanged. Finally, Figure 1(d) illustrates the convergence of the funnel method to the optimal solution.

In addition to switching condition and step types, another interesting connection exists between filter and funnel methods [28]. The filter maintains an iteratively updated list of pairs ($f$,$\theta$) to ensure progressive improvement towards optimality and feasibility, and requires heuristic tuning of two key parameters, $\lambda_\theta$ and $\lambda_f$, which define the necessary filter acceptance criteria [9]. These parameters can be problem-dependent and require careful adjustment. In contrast, we can interpret the funnel as a filter with a single entry (i.e., funnel width), which is updated only during $\theta - type$ iterations. The funnel width is initialised based on the initial infeasibility measure and a problem-independent parameter, making the funnel mechanism simpler and more adaptive in practical implementation. Unlike the filter method, the funnel does not require maintaining a list of non-dominated ($f$,$\theta$) pairs or tuning sensitive parameters like $\lambda_\theta$ and $\lambda_f$, highlighting a key practical advantages of the funnel mechanism. The similar nature of the filter and the funnel (i.e., switching condition, step types and interpretation of the funnel as a filter with single entry) help us to derive global convergence proof of the TR funnel method using the TR filter theory.

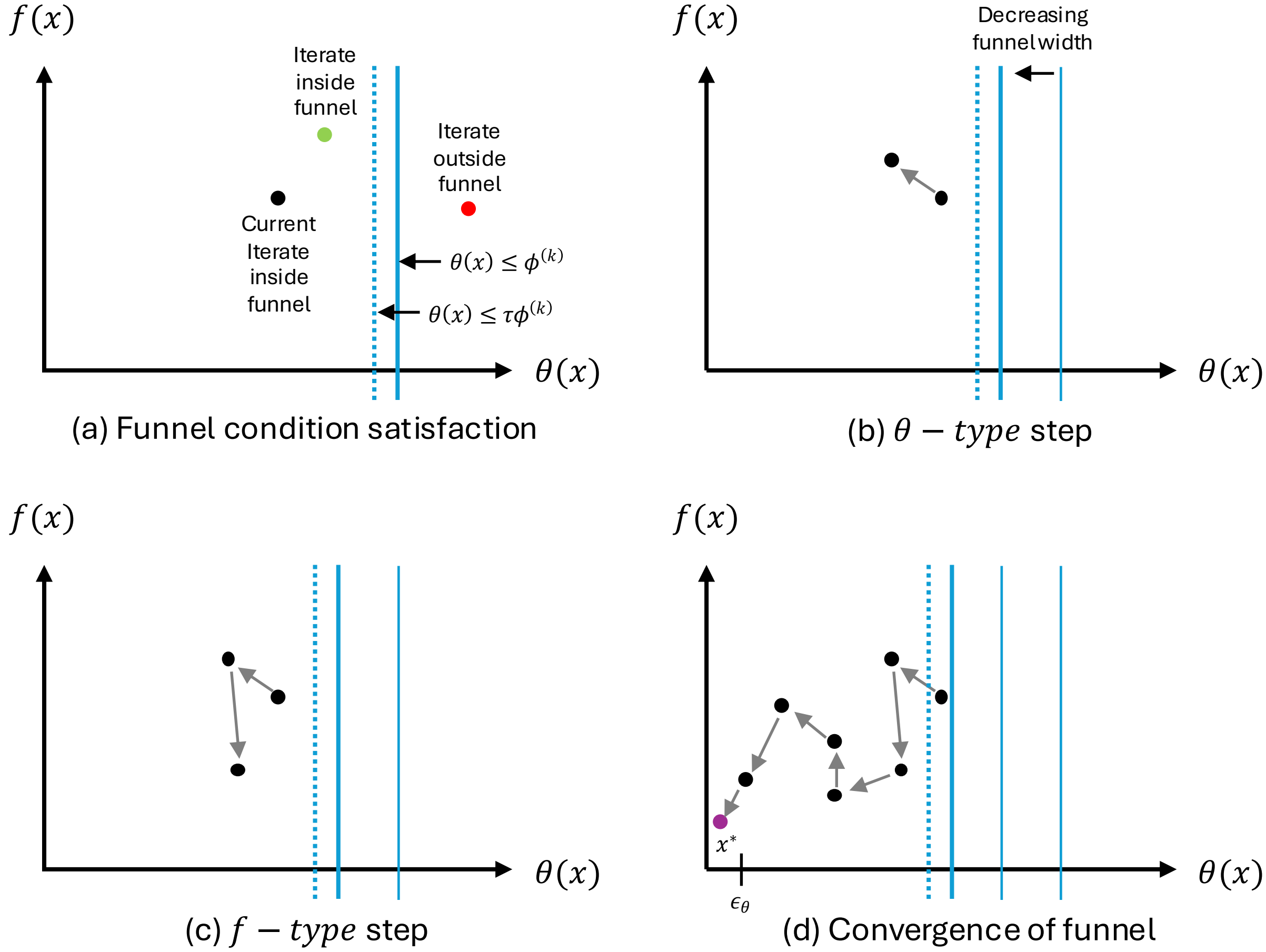

***Figure 1.*** *Funnel mechanism (solid blue line indicates necessary funnel condition (12), dashed line is funnel sufficient decrease condition (16), and light blue lines represent previous funnel widths). (a) Iterates inside the funnel (black and green) satisfy both funnel (12) and funnel sufficient decrease (16) conditions, while points outside (red) violate them. (b) When the switching condition (14) is violated but funnel (12) and funnel sufficient decrease in infeasibility (16) are satisfied (i.e., $\theta - type$ step), the funnel width is reduced to tighten feasibility requirements. (c) Both the switching (14) and Armijo decrease (15) conditions are satisfied (i.e., $f - type$ step), improving the objective function, the funnel width remains unchanged. (d) The sequence of accepted iterates progresses toward the optimal point $x^*$, as the funnel width approaches zero, demonstrating the funnel's role in ensuring convergence to feasibility and local optimality.*

## 2.5 Restoration Phase

When the compatibility check for a candidate TRSP fails, a restoration phase is activated to recover feasibility and RM consistency. Starting from the previously accepted iterate $x^{(k)}$ and current $\Delta^{(k)}$, the compatibility problem (6) is resolved iteratively by adjusting $\Delta^{(k)}$ based on the reduction ratio test (18) and re-evaluating the RM. This process continues until a compatible TRSP $(x^{(k+1)}, \Delta^{(k+1)}, r^{(k+1)}, \sigma^{(k+1)})$ is found that satisfies the funnel condition (12).

## 2.6 Algorithm

The TR funnel algorithm (as illustrated in Figure 2) is stated next.

1. Initialisation: Set $\Delta_{min}> 0$, $0 < \gamma_c < 1 < \gamma_e$, $0 < \eta_1 \leq \eta_2 < 1$, $\mu \in (0,1)$, $\gamma_s > 1/(1+\mu)$, $\kappa_\Delta \in (0,1)$, $\kappa_\phi > 1$, $\delta \in (0,1)$, $\tau \in (0,1)$, $\kappa_f \in (0,1)$, $\phi_{min} > 0$, $\eta \in (0,1)$, $\kappa_\mu > 1$, $\xi > 0$ and $\Psi \in (0,1)$.

   Set fixed suitable termination tolerances (for feasibility $\epsilon_\theta$, for criticality $\epsilon_\chi$, for compatibility $\epsilon_{comp}$, and for sampling region size $\epsilon_\Delta \geq \Delta_{min}$).

   Set initial iterate $x^{(0)}$ and other (iteratively updated) tuning parameters $\Delta^{(0)}> 0$, $\sigma^{(0)} > 0$.

   Evaluate $t(w^{(0)})$, $\theta(x^{(0)})$ and $f(x^{(0)})$. Initialise $\phi^{(0)} := max\{\phi_{min}, \kappa_\phi . \theta(x^{(0)})\}$ and $k = 0$.

2. Construct RM $r^{(k)}(w^{(k)})$ that is *κ-fully linear* in the sampling region with radius $\sigma^{(k)}$ centred at current/initial iterate $x^{(k)}$.

3. Termination and criticality check: Calculate $\chi^{(k)}$ using criticality check problem (5).

   a. A first order critical point is found as per given tolerances: if $\theta^{(k)} \leq \epsilon_\theta$, $\chi^{(k)} \leq \epsilon_\chi$, and $\sigma^{(k)} \leq \epsilon_\Delta$. STOP.

   b. A feasible point is found as per given tolerances but the progress towards optimality is slow: if $\theta^{(k)} \leq \epsilon_\theta$, $\theta^{(k-1)} \leq \epsilon_\theta$, $\Delta^{(k)}\leq \Delta_{min}$, and $\Delta^{(k-1)}\leq \Delta_{min}$. STOP.

   c. Else if criticality test (6) holds, perform criticality update (7).

4. Compatibility check: Solve the compatibility check problem (8), if $\left\| y - r^{(k)}(w) \right\| \leq \epsilon_{comp}$, go to next step. Else, call restoration phase iteration (section 2.5) to restore compatibility, set $k = k + 1$ and go to step 2. This restoration loop (steps 2-4) continues until a compatible set ($x^{(k+1)}$, $\Delta^{(k+1)}$, $r^{(k+1)}$, $\sigma^{(k+1)}$) is found.

5. Initialise the TRSP (4) using the compatible solution and then solve it to compute the step $s^{(k)}$.

6. Evaluate $\theta^{(k)}$ and $f^{(k)}$.

7. If the step satisfies funnel condition (12), continue to Step 9. Otherwise, go to step 12.

8. If switching condition (14) holds, go to step 10. Otherwise, go to step 11.

9. If Armijo-type sufficient decrease condition (15) is satisfied, the step is $\boldsymbol{f - type}$, set $x^{(k+1)} = x^{(k)} + s^{(k)}$, $\Delta^{(k+1)}= \max\{\gamma_e \left\| s^{(k)} \right\|, \Delta^{(k)}\}$, $\sigma^{(k+1)} = \sigma^{(k)}$, $\theta^{(k+1)} = \theta(x^{(k)} + s^{(k)})$, set $k = k + 1$ and go to step 2. Otherwise, go to step 12.

10. If funnel sufficient decrease condition (16) is satisfied, the step is $\boldsymbol{\theta - type}$, set $x^{(k+1)} = x^{(k)} + s^{(k)}$, $\sigma^{(k+1)} = \min[\sigma^{(k)}, \Psi\Delta^{(k)}]$, $\theta^{(k+1)} = \theta(x^{(k)} + s^{(k)})$, adjust TR size using rule (18) and Eq. (19), set $k = k + 1$ and go to step 2. Otherwise, go to step 12.

11. **Rejected step**: Set $x^{(k+1)} = x^{(k)}$, $\Delta^{(k+1)}= \gamma_c \left\| s^{(k)} \right\|$, $\theta^{(k+1)} = \theta^{(k)}$, $\sigma^{(k+1)} = \min\{\sigma^{(k)}, \Psi\Delta^{(k)}\}$, set $k = k + 1$ and go to step 2.

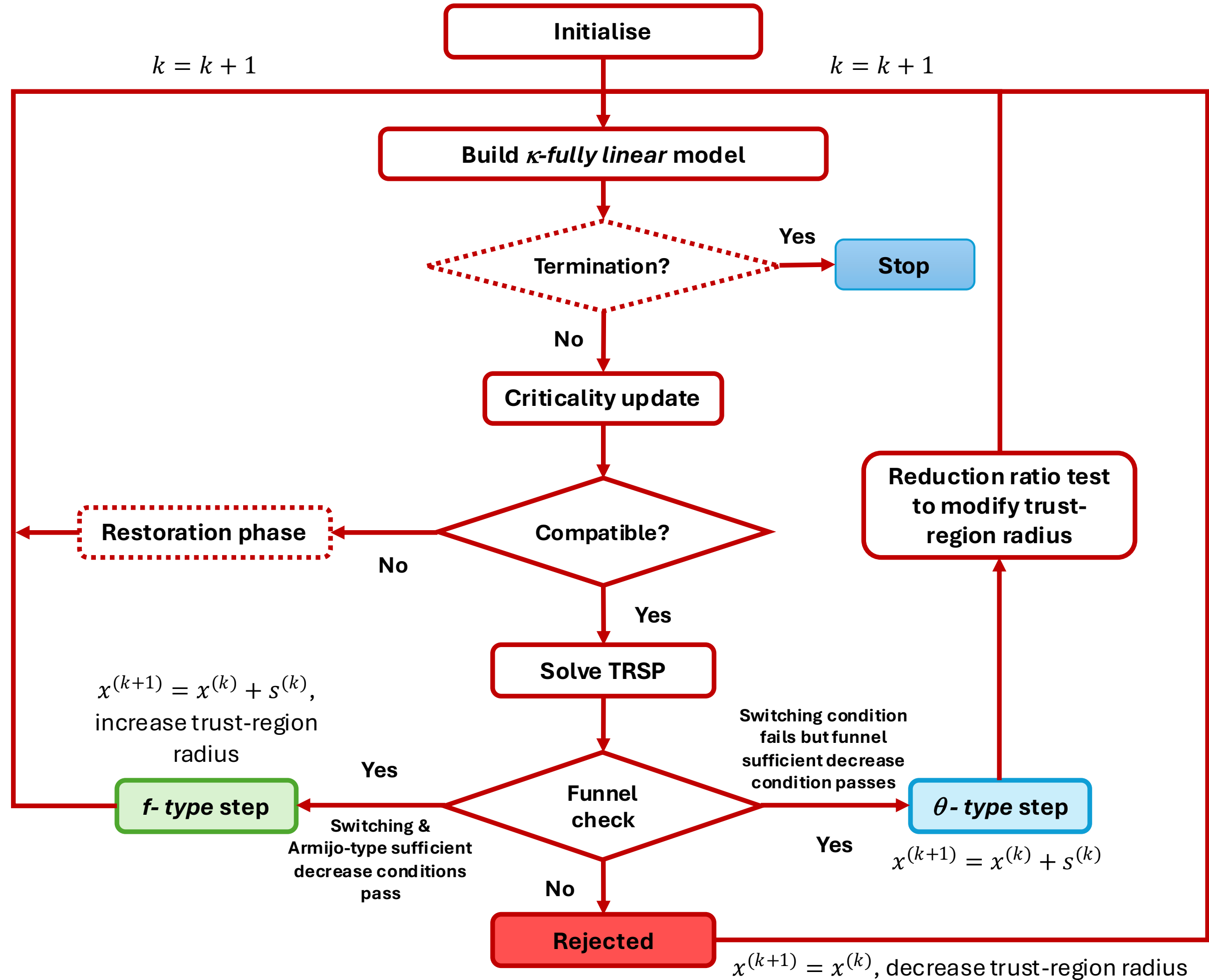

***Figure 2.*** *The trust-region funnel method for constrained grey-box optimisation. The diagram illustrates the main steps including surrogate construction, criticality and compatibility checks, trust-region subproblem (TRSP) solution, funnel acceptance and trust-region update mechanism. The iterative loops involving restoration,* $f - type$*, and* $\theta - type$ *steps guide the solution toward feasibility and local optimality of the original grey-box problem (2).*

## 3. Global Convergence

In this work, we establish convergence properties of the TR funnel algorithm. We show that the sequence of iterates generated by the TR funnel algorithm either converges to a first-order KKT point of the grey-box problem (2), given that a constraint qualification holds, or terminates unsuccessfully at a stationary point during the restoration phase after exhausting the computational budget. The proof relies on standard NLP assumptions: smoothness of the underlying functions, boundedness of level sets, and regularity of limit points [1][28]. We also assume that the TRSP solver returns a local minimiser that produces sufficient decrease relative to a feasible starting point [1][27]. Unlike classical TR convergence results, the contraction of the sampling radius $\sigma^{(k)}$ rather than the TR radius $\Delta^{(k)}$, is the primary mechanism driving convergence. The TR size $\Delta^{(k)}$ may remain bounded away from zero to allow large steps for glass-box variables (controlling globalisation). When adopting these assumptions, the RM-based constraints (constructed to be $\kappa$-*fully linear* within the sampling region) can be interpreted as equality constraints of a nonlinear subproblem, allowing the extension of classical NLP convergence proofs to the grey-box setting. We additionally assume that $\sigma^{(k)}$ is updated according to criticality test so that $\sigma^{(k)} \to 0$ as $k \to \infty$ whenever $\chi^{(k)} \to 0$. Under these conditions, it follows that there exists a subsequence $\{k_l\}$ such that $x^{(k_l)} \to x^*$, where $x^*$ denotes a first-order KKT point of problem (2). We show that $\theta^{(k_l)}$, $\chi^{(k_l)}$, and $\sigma^{(k_l)}$ all converge to zero, implying that $x^*$ satisfies first-order optimality conditions. The detailed convergence proof is stated in Appendix 1.

## 4. Implementation

The TR funnel method is implemented in Pyomo, an open-source Python framework for mathematical programming [34]. Pyomo enables integration of the TR funnel framework with a range of NLP solvers; by default, subproblems are solved using IPOPT 3.14.13 with MA27 [7]. Derivatives for the criticality check are computed via the AMPL pseudo-solver *gjh*.

Numerical experiments were conducted on a Windows 11 system with an 11th Gen Intel(R) Core(TM) i7-1165G7 @ 2.80 GHz (8 CPUs, ~1.7 GHz) and 16 GB RAM. No explicit parallelisation or multithreading is employed beyond the default solver behaviour. The complete Pyomo code, including examples, tuning parameters, and solver settings, is available in our public GitHub repository [35].

## 5. Numerical Results

We evaluate the new TR funnel method against the TR filter algorithm [9] across a set of constrained nonlinear grey-box test problems, engineering design benchmarks, and process models. The benchmark set comprises the Colville and Himmelblau functions, the Loeppky problem, the Wing Weight and Welded Beam structural design problems, and two process optimisation case studies: Williams–Otto process and biomass-based hydrogen production flowsheet problems. Each NLP is reformulated to include black-box components $t_i(w)$ that mimic simulation-based outputs. Both funnel and filter methods used identical initial points, parameter settings where applicable, and subproblem solvers. The method-to-method performance is compared in terms of the number of iterations, black-box evaluation count, and computational time. The iteration-wise convergence behaviour is also analysed for selected problems using metrics such as the infeasibility measure, criticality measure, sampling region size, and step size.

### 5.1 Numerical Benchmarks

We first assess the TR funnel algorithm on three numerical constrained grey-box optimisation problems, which are adapted from classical analytical test functions: Colville, Himmelblau, and Loeppky functions. Each algorithm is tested using five RM forms (linear, quadratic, simple quadratic excluding variable interaction terms, Gaussian process, and Taylor series adapted from [26]) to investigate the influence of RM accuracy on algorithmic performance. The modified Colville problem has dimensions $(n_w, n_y, n_z) = (4,4,1)$, yielding a total of nine variables and six nonlinear inequality constraints. The modified Himmelblau problem has $(n_w, n_y, n_z) = (3,2,5)$ variables, leading to ten variables and three nonlinear equality relations coupled with boundary conditions. The Loeppky problem is the smallest with $(n_w, n_y, n_z) = (3,1,4)$ variables but still exhibits strong nonlinearity and variable interaction. Mathematical grey-box formulations are provided in the supplementary material (Annexe A). These formulations allow experimentation with varying dimensionality and constraint structure while retaining the features of practical engineering problems.

All tests began with identical initial objective and infeasibility values for each problem type. Both funnel and filter algorithms consistently converged to the same optimal values across all problems and RM forms, indicating that the RM choice did not impact the final accuracy in these cases.

Table 1 reports the final objective value, number of external function evaluations, and total CPU time in seconds (wall-clock runtime measured using Python) for each RM-algorithm configuration. In most cases, the TR funnel algorithm required fewer evaluations and less CPU time than the TR filter for the same RM forms. The advantage was most pronounced for Colville with quadratic models and Himmelblau with linear and standard quadratic RMs, where the TR funnel variant achieved convergence with up to 4-135 times fewer evaluations than the TR filter. In contrast, for highly accurate Gaussian process and Taylor series RMs, both methods exhibited almost identical costs, with the TR funnel showing marginal time savings due to slightly fewer external evaluations. For Loeppky, both methods converged requiring very few evaluations, with no performance difference. These patterns align with convergence theory: when RM accuracy is high, both TR methods behave similarly; when RM accuracy is weak, the funnel's feasibility-handling strategy limits excessive feasibility correction cycles, reducing external evaluation overhead. The reduced computational cost reflects a key element of the TR funnel convergence proof: once the feasibility tube contracts sufficiently to contain the iterates, feasibility is retained, enabling more optimality steps.

***Table 1.*** *Performance metrics for the proposed TR funnel and the classical TR filter algorithms on numerical benchmark problems (Colville, Himmelblau, and Loeppky). Each problem was solved using five reduced RM forms: linear (L), quadratic (Q), simple quadratic (SQ), Gaussian process (GP), and Taylor series (TS). Both algorithms converged to identical objective values, demonstrating surrogate-independence of the final solution. The TR funnel method generally required fewer external evaluations and less CPU time, particularly for Colville (with Q and SQ RMs) and Himmelblau (with all RM forms), indicating improved computational efficiency while maintaining comparable accuracy to TR filter.*

| **Problem** | **Reduced Model Form** | **TR Funnel** | | | **TR Filter** | | |
|---|---|---|---|---|---|---|---|
| | | **Final Objective** | **External Evaluations** | **CPU Time (s)** | **Final Objective** | **External Evaluations** | **CPU Time (s)** |
| Colville | L | 10122.49 | 50412 | 1523.29 | 10122.49 | 36724 | 1125.10 |
| | Q | 10122.49 | 3012 | 32.91 | 10122.49 | 407988 | 4183.18 |
| | SQ | 10122.49 | 1984 | 30.30 | 10122.49 | 255724 | 4108.85 |
| | GP | 10122.49 | 204 | 13.32 | 10122.49 | 204 | 13.63 |
| | TS | 10122.49 | 14 | 2.22 | 10122.49 | 14 | 2.13 |
| Himmelblau | L | -25822.94 | 1671 | 97.96 | -25822.95 | 6933 | 346.78 |
| | Q | -25822.94 | 1101 | 30.22 | -25822.95 | 8649 | 166.74 |
| | SQ | -25822.95 | 603 | 15.15 | -25822.95 | 603 | 13.62 |
| | GP | -25822.95 | 147 | 10.29 | -25822.95 | 243 | 19.98 |
| | TS | -25822.95 | 13 | 3.37 | -25822.95 | 14 | 4.59 |
| Loeppky | L | -1.55E-07 | 13 | 1.59 | -1.55E-07 | 13 | 1.53 |
| | Q | -1.55E-07 | 34 | 1.57 | -1.55E-07 | 34 | 1.68 |
| | SQ | -1.55E-07 | 25 | 1.69 | -1.55E-07 | 25 | 1.64 |

| | GP | -1.55E-07 | 25 | 2.83 | -1.55E-07 | 25 | 2.98 |
|---|---|---|---|---|---|---|---|
| | TS | -1.55E-07 | 7 | 1.75 | -1.55E-07 | 7 | 1.82 |

To better understand iteration-level behaviour, Table 2 reports the step-type distributions. For the TR funnel, 96.6% of all steps across problems were accepted (97.6% $f - type$, remainder $\theta - type$), with just 2.5% rejected and 0.9% restoration steps. The distribution varied with problem geometry:

- Colville: Very high acceptance (98.9%), dominated by $f - type$ steps (99.5%), indicating rapid progress toward optimality with minimal feasibility corrections.
- Himmelblau: Lower acceptance (72%) and a higher rejection rate (26.1%), reflecting the tighter feasible region and more challenging nonlinear geometry.
- Loeppky: 100% acceptance with a two-thirds $f - type$ split, consistent with its near-linear feasible region.

The TR filter's acceptance rates were even higher overall (99.7%) but with a slightly lower $f - type$ proportion (95.9%). Notably, the Himmelblau problem under the TR filter required far more $\theta - type$ steps (69% as compared to 27% with TR funnel), inflating iteration counts and evaluations despite high acceptance. This is consistent with the TR filter's convergence proof, which predicts that weaker RMs increase the $\theta - type$ burden, whereas the funnel's shrinking feasibility tube mitigates this effect.

***Table 2.*** *Step-type distribution for the proposed TR funnel and the classical TR filter algorithms applied to numerical benchmark problems (Colville, Himmelblau, and Loeppky). Each problem was solved using five RM forms: linear (L), quadratic (Q), simple quadratic (SQ), Gaussian process (GP), and Taylor series (TS). The TR funnel algorithm consistently exhibits a higher share (combined for all problems) of* $\boldsymbol{f - type}$ *steps and fewer* $\boldsymbol{\theta - type}$ *steps, reflecting a faster smoother convergence behaviour. In contrast, the TR filter shows an increased proportion of* $\boldsymbol{\theta - type}$ *steps, particularly for low-fidelity surrogates (L and Q), consistent with its filter-based feasibility correction mechanism. Across all problems, both algorithms converged reliably, with TR funnel showing comparable and often improved computational efficiency.*

| **Problem** | **Reduced Model Form** | **TR Funnel** | | | | **TR Filter** | | | |
|---|---|---|---|---|---|---|---|---|---|
| | | **f-type** | **θ-type** | **Rejected** | **Restoration** | **f-type** | **θ-type** | **Rejected** | **Restoration** |
| Colville | L | 2502 | 7 | 3 | 10 | 1825 | 6 | 0 | 4 |
| | Q | 41 | 4 | 0 | 2 | 6355 | 4 | 2 | 22 |
| | SQ | 42 | 1 | 1 | 6 | 6383 | 7 | 1 | 2 |
| | GP | 1 | 1 | 1 | 2 | 0 | 3 | 0 | 2 |
| | TS | 1 | 0 | 1 | 2 | 1 | 1 | 0 | 2 |
| Himmelblau | L | 109 | 15 | 42 | 1 | 79 | 491 | 11 | 1 |
| | Q | 26 | 3 | 22 | 1 | 186 | 75 | 0 | 1 |
| | SQ | 0 | 24 | 0 | 1 | 0 | 24 | 0 | 1 |
| | GP | 0 | 4 | 1 | 1 | 3 | 5 | 1 | 1 |

| | TS | 0 | 4 | 2 | 1 | 2 | 5 | 0 | 1 |
|---|---|---|---|---|---|---|---|---|---|
| Loeppky | L | 2 | 1 | 0 | 0 | 1 | 2 | 0 | 0 |
| | Q | 2 | 1 | 0 | 0 | 1 | 2 | 0 | 0 |
| | SQ | 2 | 1 | 0 | 0 | 1 | 2 | 0 | 0 |
| | GP | 2 | 1 | 0 | 0 | 1 | 2 | 0 | 0 |
| | TS | 2 | 1 | 0 | 0 | 1 | 2 | 0 | 0 |

High-fidelity (Gaussian process and Taylor series) RMs out-performed lower-fidelity polynomial (linear, standard quadratic and simple quadratic) RMs in all performance metrics. For example, the Colville problem using Gaussian process RMs in the TR funnel implementation required only 2 accepted steps plus 1 rejected and 1 restoration step. In contrast, linear RMs could require thousands of steps (e.g., the Colville problem required more than 1800 $f - type$ steps) due to repeated model refinement and feasibility corrections. Overall, accurate and high-fidelity local models yielded high TR acceptance rates and minimal feasibility restoration step usage.

## 5.2 Engineering Design Benchmarks

We now evaluate the TR funnel algorithm on two engineering design problems: the Wing Weight and Welded Beam problems, both re-formulated as constrained grey-box optimisation problems (Annexe A). These benchmarks are representative of practical engineering applications where design variables show nonlinear interactions. We also compare the TR funnel results with the TR filter method for completeness, but our primary focus remains on assessing TR funnel's performance characteristics, particularly in the context of its convergence behaviour.

**Wing Weight Optimisation**

The Wing Weight problem models the structural design of a light aircraft wing, where the objective is to minimise total weight while satisfying geometric and material constraints. In our modified formulation, the problem has $(n_w, n_y, n_z) = (2,1,8)$ variables, with an initial objective value of 251.85 and an initial infeasibility measure of 7.75, reflecting a highly infeasible starting point.

Using the TR funnel algorithm, all RM forms converged to the same local optimum with a final objective of 123.25 in 15 iterations (Figure 3). For this problem, RM choice does not affect convergence behaviour, although it influences computational cost. Taylor series RM required the fewest external evaluations (32), whereas the standard quadratic approximation required the most (106). Gaussian process RM achieved the solution with similar iteration counts but incurred the highest CPU time (~15 seconds).

Across all RM forms, the algorithm executed two initial restoration steps, which reduced infeasibility sharply from 7.75 to 2.55 and then 1.32. After this, the funnel accepted consecutive $f - type$ steps that reduced the objective while maintaining feasibility. The infeasibility measure remained below $10^{-8}$ after the second step and its average across RM forms is shown in Figure 3(a). Convergence to the local optimum was achieved after one rejected step. Results highlight the two-phase convergence typical of funnel methods, rapid infeasibility reduction to constrain the infeasibility within the funnel, followed by monotonic objective decrease and non-increasing funnel width under controlled TR updates.

The TR filter method with all RM forms also converged to the same objective value in 15 iterations as shown in Figure 3(b). The different step pattern (a $\theta - type$ step after two initial restoration steps, followed by a sequence of $f - type$ steps until convergence) resulted in a slightly more external evaluations when using quadratic RM (110 compared to 106 with TR funnel method).

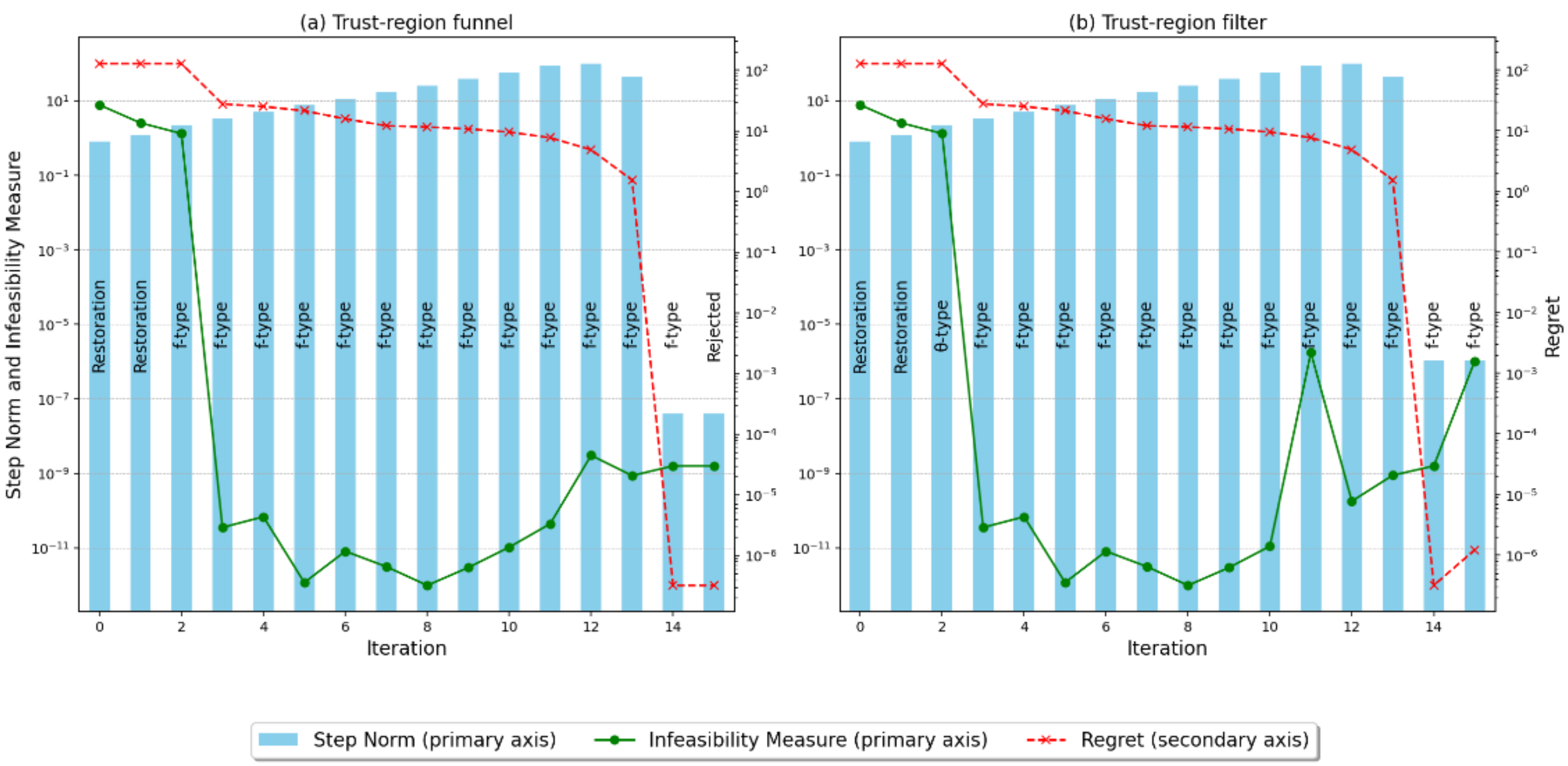

***Figure 3.*** *Iteration-wise convergence behaviour for the Wing Weight optimisation problem using the (a) TR funnel and (b) TR filter algorithms. Blue bars represent the step norm, green solid lines the infeasibility measure, and red dashed lines the regret (absolute deviation from the true objective). All RM forms showed same convergence behaviour for each algorithm. Both methods rapidly restore feasibility within the first few iterations and then proceed predominantly with $f - type$ steps. The TR funnel shows smoother feasibility reduction, whereas the TR filter exhibits similar trends but with marginally higher oscillations in infeasibility near termination.*

**Welded Beam Optimisation**

The Welded Beam problem represents a structural design scenario in which the objective is to minimise fabrication cost subject to constraints on shear stress, normal stress, deflection, and buckling load. The modified formulation has $(n_w, n_y, n_z) = (4,1,0)$ variables, with an initial objective value of 1 and an initially high infeasibility measure of 9.1.

Figure 4 shows the iteration-wise convergence metrics (step norm, infeasibility measure, regret, and criticality measure) for the Welded Beam optimisation. Using the TR funnel algorithm, all RM forms converged to the same local optimum with a final objective of 1.72. Taylor series and Gaussian process RM converged in only 7 iterations (with same convergence behaviour), while quadratic RM required slightly more, with the linear RM being slowest, requiring 13 iterations. The standard quadratic RM required the most external evaluations (123 in 8 iterations), whereas the Taylor series RM was most efficient, requiring only 16.

The step-type sequence again illustrates the theoretical mechanism of the funnel method. Early restoration and $\theta - type$ steps rapidly reduced infeasibility from 9.1 to below $10^{-1}$, after which the method switched to $f - type$ steps, steadily reducing the objective while keeping feasibility safeguarded. Convergence was finalised by a rejected, or $\theta - type$, step.

In contrast, the TR filter algorithm displayed a different convergence profile, especially with polynomial RM forms. With linear RM, it required 50 iterations and 256 external evaluations to reach a local solution. Even with quadratic RM, the TR filter method required 139 calls in 9 iterations, more than the TR funnel method under equivalent reduced models. This difference supports the claim from the funnel convergence proof that the explicit separation of feasibility and optimality mechanisms in the funnel enables efficient feasibility restoration when starting from infeasible points.

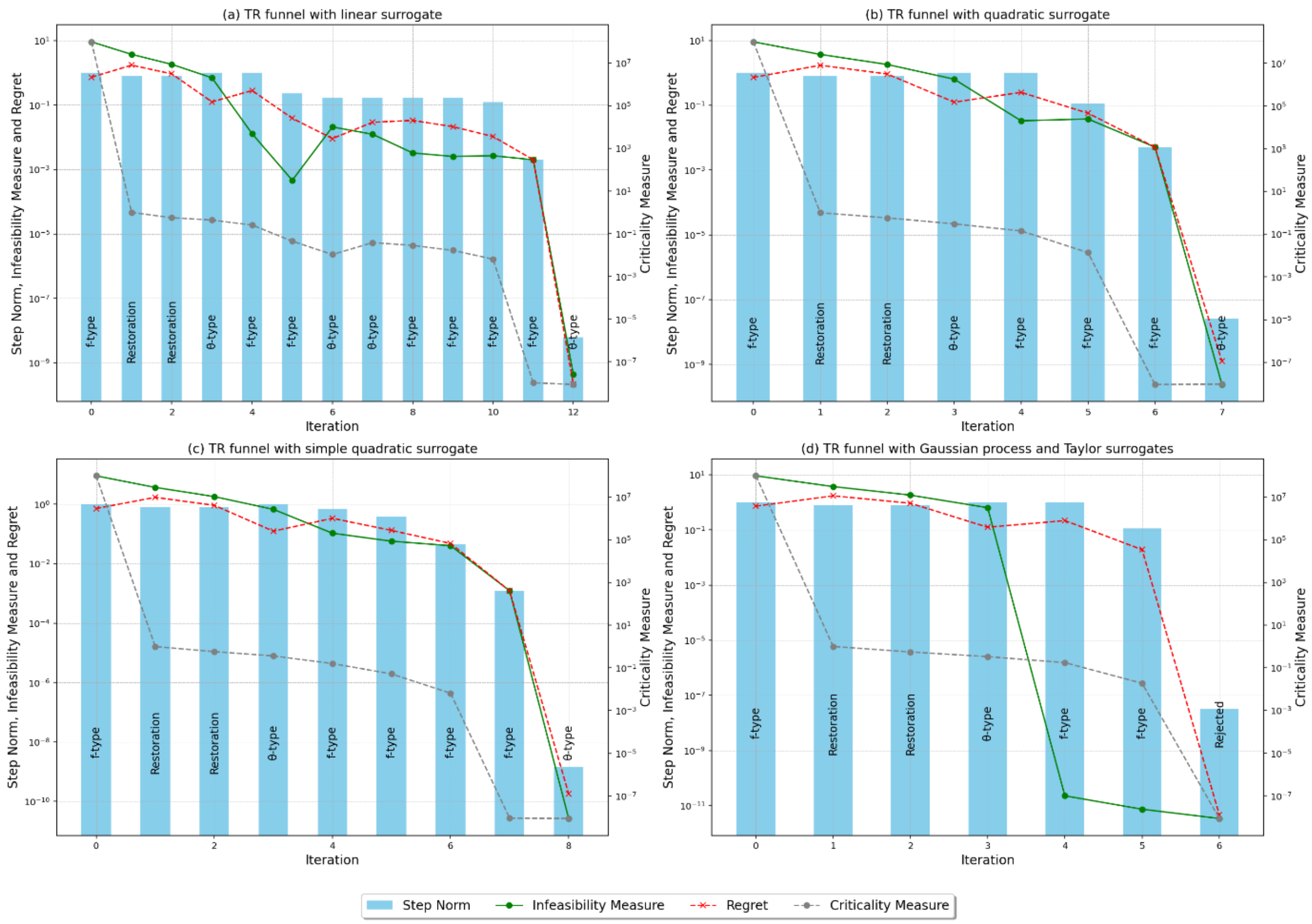

***Figure 4.*** *Iteration-wise convergence of the Welded Beam optimisation problem using the TR funnel algorithm with different RM forms: (a) linear, (b) quadratic, (c) simple quadratic, and (d) Gaussian process and Taylor series. Blue bars represent step norms, while green solid, red dashed, and grey dashed lines denote infeasibility measure, regret, and criticality measure, respectively. All RM forms achieve rapid feasibility restoration followed by monotonic decrease in criticality measure, converging to the same local optimum (regret) within 6–13 iterations. The Taylor and Gaussian Process RMs converge fastest, while low-order polynomial models require more iterations for comparable accuracy.*

## 5.3 Process Optimisation Case Studies

In this section, we demonstrate our TR funnel method on two process systems case studies: the Williams-Otto process and biomass-based hydrogen production via gasification.

**Williams-Otto Process**

The Williams-Otto (Figure 5) is a well-known benchmark problem in process optimisation [36]. Two feed streams, $F_A$ and $F_B$, are combined with a recycle stream $F_R$ and introduced into a continuous stirred-tank reactor (CSTR). Three sequential second-order reactions occur:

$$A + B \xrightarrow{k_1} C,$$

$$B + C \xrightarrow{k_2} P + E,$$

$$P + C \xrightarrow{k_3} G,$$

where $P$ is the desired product, $C$ is an intermediate, $E$ is a by-product and $G$ is a waste component. The reactor kinetics are treated as a black box mapping $t(w)\colon \mathbb{R}^6 \rightarrow \mathbb{R}^3$.

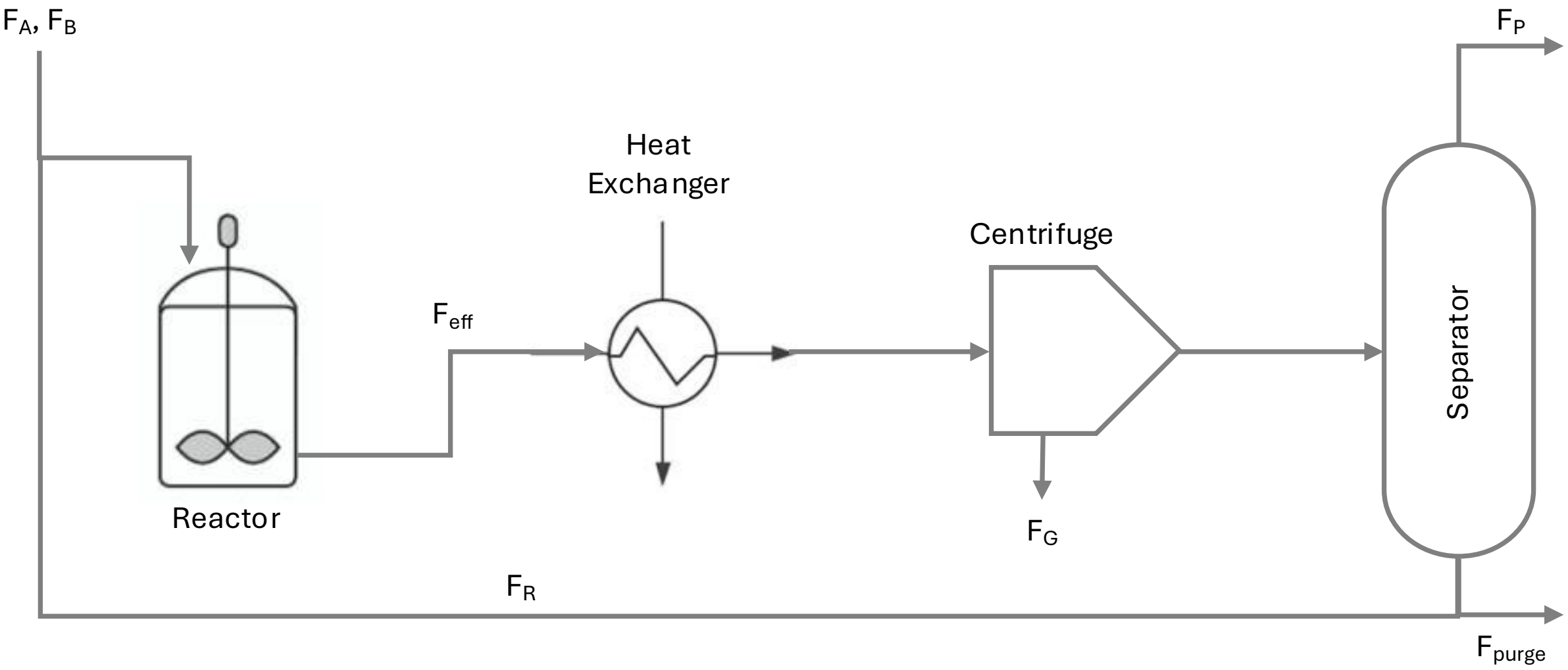

***Figure 5.*** *Process flowsheet of the Williams–Otto reaction system, consisting of a CSTR, a shell-and-tube heat exchanger, a centrifuge (decanter), and a distillation separator. Feed streams $F_A$ and $F_B$ enter the reactor, where three consecutive reactions produce intermediates and products. The effluent $F_{eff}$ is cooled, separated to remove waste gas $F_G$, and split into product ($F_P$) and purge ($F_{purge}$) streams, with recycle $F_R$ returning to the reactor.*

The effluent stream $F_{eff}$ is cooled in a heat exchanger and then sent to a centrifuge to separate component $G$ in the stream $F_G$. The remaining mixture is further processed to recover product $P$ in the overhead stream $F_P$. Due to the presence of an azeotrope, 10 wt % of component $E$ remains in the bottoms, which is split into purge stream $F_{purge}$ and recycle stream $F_R$. The resulting recycle stream is mixed with fresh feed and re-introduced into the reactor. The detailed grey-box process optimisation model with $(n_w, n_y, n_z) = (6,3,21)$ variables, is provided in Annexe A (AN6).

Starting from an initial objective value of $-11.54$ with an infeasibility measure of 6.1, the TR funnel method with Taylor series RM converges to a locally optimal solution of $-121.03$ in 34 iterations requiring 195 external evaluations. Due to feasible initial values and the local accuracy of Taylor model, the optimisation proceeds with a sequence of $f - type$ steps before converging after several rejected steps (Figure 6(a)). For comparison, the same problem solved using the TR filter method (Figure 6(b)) begins with two initial $\theta - type$ steps, followed by consecutive $f - type$ steps. After a final $\theta - type$ step, the TR filter method reaches a similar local solution in 37 iterations using 222 external black-box evaluations.

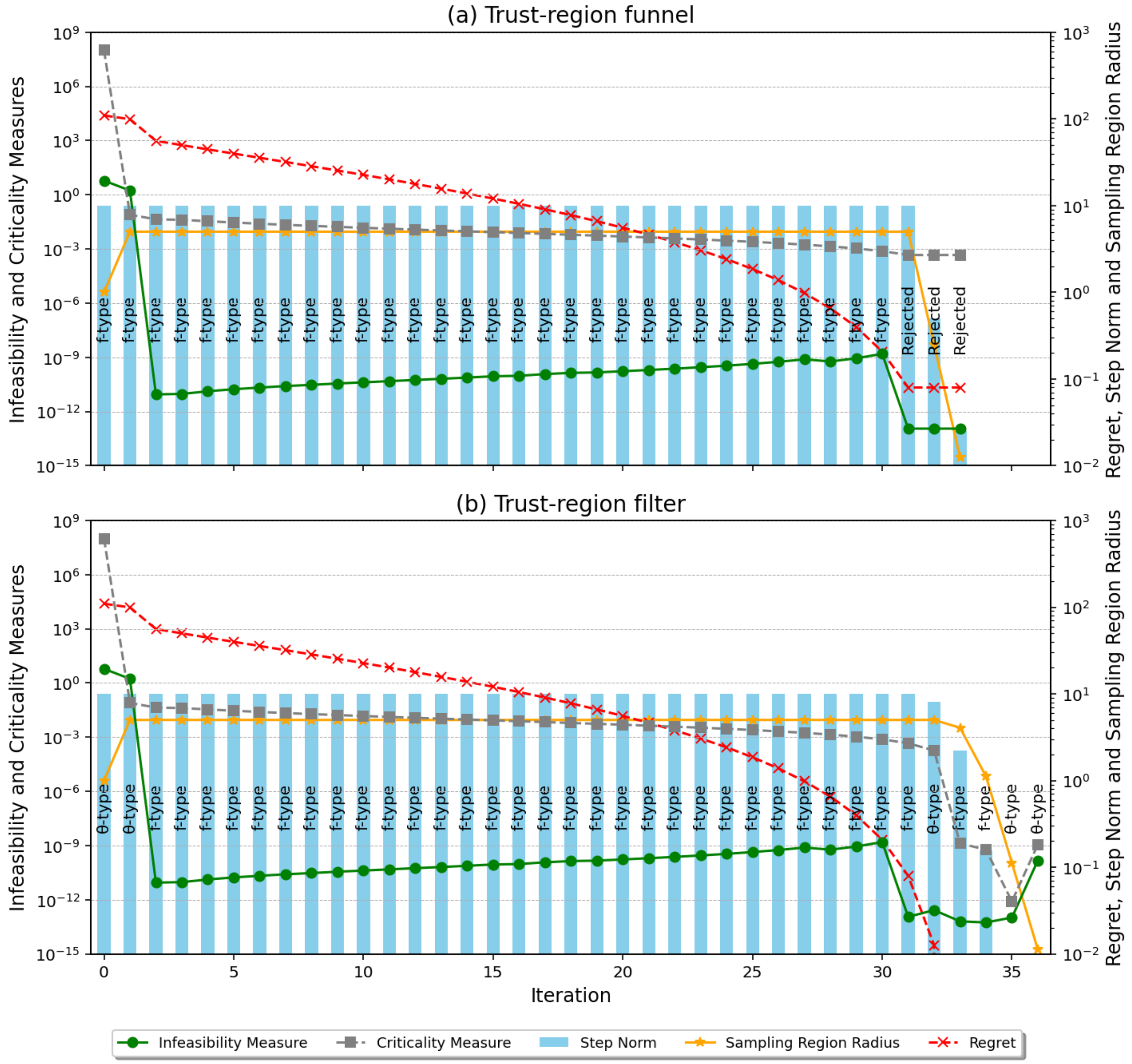

***Figure 6.*** *Iteration-wise convergence behaviour for the Williams–Otto process optimisation using the (a) TR funnel and (b) TR filter algorithms. Blue bars represent step norms, green solid lines the infeasibility measure, red dashed lines the regret (deviation from true objective), yellow solid lines the sampling region radius, and grey dashed lines the criticality measure. Both methods achieve rapid feasibility recovery followed by steady objective improvement, converging to the same local optimum. The TR funnel exhibits smoother infeasibility decay, faster convergence and slightly fewer* $\theta - type$ *steps compared to the TR filter, demonstrating its improved performance.*

**Biomass-based Hydrogen Production**

We consider a process systems optimisation problem for minimising the capital (CAPEX) and operating (OPEX) expenditures of a biomass gasification hydrogen production plant. The process flowsheet (Figure 7) consists of a gasification unit, a combined heat and power (CHP) engine, a syngas cleaner, a water–gas shift (WGS) reactor, and a hydrogen separation unit.

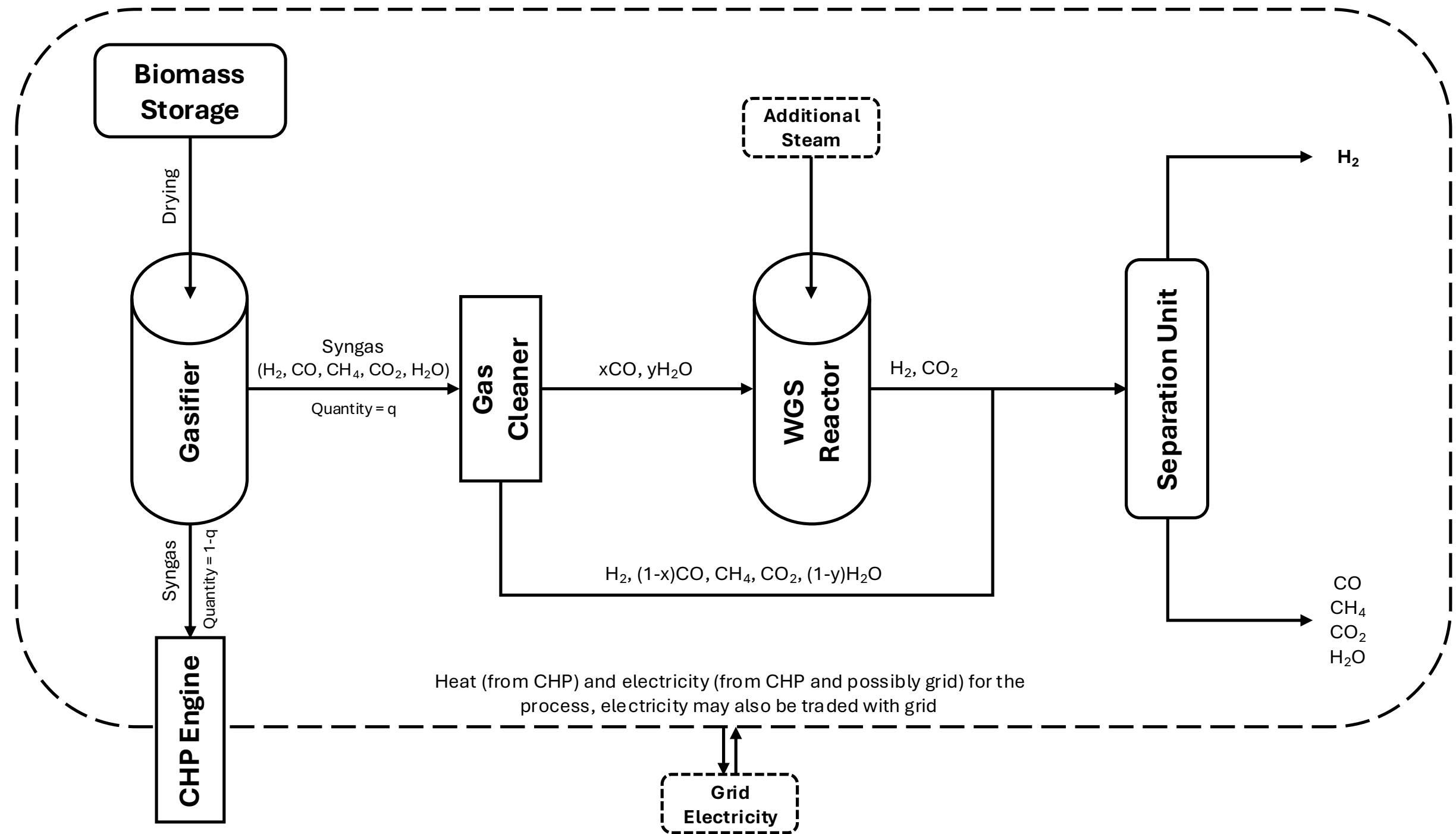

***Figure 7.*** *Flowsheet of the biomass-based hydrogen production process, comprising a gasifier, gas cleaner, WGS reactor, and separation unit. Biomass is converted to syngas (*$H_2$*,*$CO$*,* $CH_4$*,* $CO_2$*,* $H_2O$*) in the gasifier, partially processed in a CHP engine, and further shifted in the WGS reactor to enhance hydrogen yield. The separation unit seaparates pure* $H_2$ *from residual gases (*$CO$*,* $CH_4$*,* $CO_2$*,* $H_2O$*). Heat and electricity from the CHP and grid are integrated for process energy needs and may also be traded externally.*

This techno-economic grey-box optimisation problem had $(n_w, n_y, n_z) = (8,4,71)$ variables and 69 constraints. The reaction kinetics of the gasifier was partially black-boxed, while the downstream WGS reactor was fully treated as a black-box, mapping $t(w): \mathbb{R}^8 \rightarrow \mathbb{R}^4$. Black-box inputs included gasifier temperature, syngas compositions, and WGS outlet conditions. The glass-box part of the model retained the mass and energy balances, economic equations, power generation, cost correlations, and other process-wide constraints.

The TR funnel method with Taylor series RM for the 4 black-box components ($n_y = 4$) successfully solved the biomass-based hydrogen production problem to local optimality in 193 iterations. The infeasibility measure, initially at 15, was reduced to $\sim 10^{-15}$ within the first two iterations and remained negligible ($<10^{-10}$) thereafter (Figure 8). After an $f - type$ and then $\theta - type$ step in the beginning, the algorithm executed $f - type$ steps consecutively, steadily reducing regret from $\sim 10^3$ to below $10^{-2}$. Convergence was achieved at an objective value of $9,816,106, coinciding with an expected sharp reduction during the step norm in last two steps. The TR filter method also converged in 193 iterations with identical step profile. Both TR Funnel and TR filter method required 1544 external evaluations to converge.

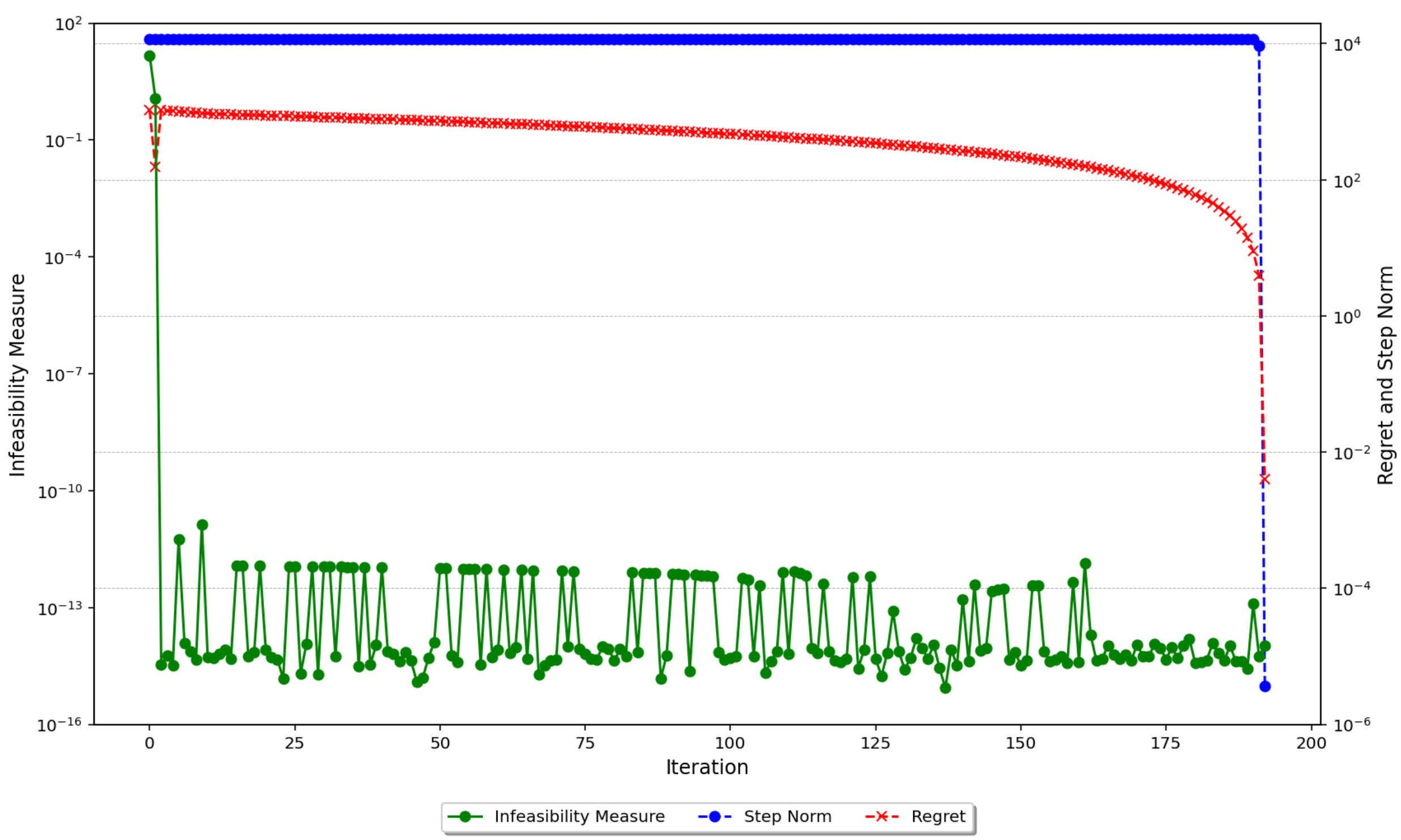

***Figure 8.*** *Iteration-wise convergence of the biomass-based hydrogen production process optimisation using the TR funnel algorithm. The green solid lines represents the infeasibility measure (left axis), while the blue and red dashed correspond to the step norm and regret (right axis), respectively. Feasibility is restored rapidly within the first few iterations, and the infeasibility remains below* $10^{-10}$ *thereafter. The regret and step norm decrease monotonically, demonstrating smooth objective improvement and convergence to a local optimum after 193 iterations. This behaviour confirms the TR funnel's robustness and consistent feasibility enforcement for large-scale grey-box systems.*

## 6. Conclusions and Recommendations

We have presented a TR funnel method for RM-based grey-box optimisation, accompanied by a rigorous convergence proof to first-order stationary point. Unlike TR filter methods which manage both objective and infeasibility measure reduction as multiple objectives requiring more tuning parameters, the TR funnel monotonically decreases the infeasibility measure while facilitating the TR globalisation strategy to achieve local optimality. Across a range of benchmark and process optimisation problems, the TR funnel method showed robust convergence and performance comparable to the classical TR filter approach. As per the convergence proof, the computational experiments demonstrate that the TR funnel approach guides the optimisation (via restoration or $\theta - type$ steps) toward the feasible space in the beginning, followed by the reliable convergence to local optima (using $f - type$ steps), finally terminating at rejected or $\theta - type$ steps (resulting in smaller steps at the end). These findings establish the TR funnel method as a theoretically sound and practically effective tool for large-scale grey-box optimisation. Future work will explore Hessian-informed methods and sampling region update mechanisms to further enhance global convergence of the TR funnel and the TR filter methods for systems with multiple black-box components.

### Data Availability Statement

Pyomo implementation of this work (algorithms and case tests) is accessible through our public GitHub repository for the TR Filter/Funnel solver [35].

Supplementary material contains detailed mathematical formulation of the performed case tests.

## Appendix A: Global Convergence of Trust-Region Funnel Algorithm

We begin by restating the definition of constraint qualification and first-order KKT conditions, followed by a modified set of standard NLP assumptions to accommodate the RM-based TRSP formulation. Finally, we establish global convergence properties of the TR funnel algorithm under mild assumptions.

### Definitions

**Definition 1** (Constraint Qualification). *Consider a feasible point* ${x'}^T = \left[{w'}^T, {y'}^T, {z'}^T\right]$ *and a constraint set* $C = \left\{x' \colon c_i(x') = 0, g_j(x') \leq 0;\ \forall i \in \mathcal{E}, j \in \mathcal{I}\right\}$ *where the equality constraint is defined as* $c(x') = (h(x')^T, (y' - t(w')^T)$ *and the active inequality constraint set is* $\mathcal{A}(x') = \left\{j \in \mathcal{I} \colon g_j(x') = 0\right\}$. *We say that the Mangasarian–Fromowitz constraint qualification (MFCQ) holds for* $C$ *at* $x'$ *if the gradients of the equality constraints* $\{\nabla c_i(x')\}_{\,i\in\mathcal{E}}$ *are linearly independent (i.e., they have a full column rank) and there exists a direction* $d \in \mathbb{R}^n$ *such that* $\nabla c_i(x')^T d = 0;\ \forall i \in \mathcal{E}$ *and* $\nabla g_j(x')^T d < 0;\ \forall\, j \in \mathcal{A}(x')$.

**Definition 2** (First-Order KKT Conditions). *If MCFQ holds at an optimal point* ${x^*}^T = \left[{w^*}^T, {y^*}^T, {z^*}^T\right]$ *of the nonlinear TRSP, then there exist Lagrange multipliers* $\lambda^* \in \mathbb{R}^m$ *and* $\mu^* \in \mathbb{R}^n$ *such that the following first-order necessary optimality conditions are satisfied:*

$$\begin{gathered}\nabla_x \mathcal{L}(x^*, \lambda^*, \mu^*) = \nabla f(x^*) + \nabla \mathrm{c}(x^*)\lambda^* + \nabla g(x^*)\mu^* = 0, \\ c(x^*) = 0, g(x^*) \leq 0, \mu^* \geq 0, g(x^*)^T \mu^* = 0,\end{gathered} \tag{A1}$$

*where* $c(x^*) = [h(x^*)^T, (y^* - t(w^*)^T)]$. *These (A1) are collectively referred to as KKT conditions.*

**Definition 3** (Limit Point). *A point* $x^*$ *is a limit point of a sequence* $\left\{x^{(k)}\right\}$ *if there exists a subsequence* $\{x^{(k_j)}\}$ *such that* $x^{(k_j)} \to x^*$.

### Assumptions

We introduce a set of standard assumptions that ensure the TR funnel algorithm is well-posed and guarantee convergence properties analogous to those in classical NLP theory.

**AS1** (Smoothness and Bounded Derivatives). All functions (objective and constraints) are twice continuously differentiable on an open set containing $X$, with uniformly bounded first and second derivatives.

**AS2** (Bounded Iterates). The sequence of iterates $\left\{x^{(k)}\right\}$ remains within a nonempty compact domain $X \subset \mathbb{R}^{m+n+p}$.

**AS3** (Reduced Model Accuracy). The reduced models $r^{(k)}(w)$ used in TRSP (4) are $\kappa$-*fully linear* at each iteration $k$, and have uniformly bounded second derivatives.

**AS4** (Constraint Qualification). MFCQ holds for the constraint set of problem (3) at all limit points of the sequence $\left\{x^{(k)}\right\}$.

**AS5** (Fraction of Cauchy Decrease). The solution $s^{(k)}$ to the TRSP (4) satisfies a fraction of Cauchy decrease condition:

$$f\left(x^{(k)}+d^{(k)}\right)-f\left(x^{(k)}+s^{(k)}\right) \geq \kappa_{tmd}\chi^{(k)} min\left[\frac{\chi^{(k)}}{\beta^{(k)}}, \Delta^{(k)}\right], \tag{A2}$$

where $\kappa_{tmd}$ is a constant and $\beta^{(k)} > 1$ is a bounded sequence.

**AS6** (Regularity). For some constants $\delta_d > 0$ and $\kappa_{usc} > 0$, independent of $k$, whenever $\theta^{(k)} \leq \delta_d$, there exists a feasible solution $x''^T = \left[w''^T, y''^T, z''^T\right]$ to the compatibility problem (8) such that $y'' - r^{(k)}(w'') = 0$ and

$$\left\|d^{(k)}\right\| \leq \kappa_{usc}\theta^{(k)}. \tag{A3}$$

It is further assumed that whatever method is used to solve (8) will always compute such a solution.

### Convergence to Feasibility

Convergence to feasibility is largely independent of the specific globalisation strategy employed (trust-region or line-search) and instead relies on three key elements of the algorithm: (i) the switching condition, (ii) the boundedness of the objective function $f$, and (iii) the funnel mechanism that enforces progressive reduction of the infeasibility measure. Theorem 1 establishes that the algorithm generates a sequence of iterates that approaches feasibility. The result is derived under the implicit assumption that the inner restoration phase of the TR funnel method terminates successfully, and the TR funnel algorithm produces an infinite sequence of iterates [28].

**Theorem 1.** *Assume that the TR funnel algorithm does not terminate finitely at a KKT point of problem (3), and consider the sequences* $\{\phi^{(k)}\}$, $\{\theta^{(k)}\}$, $\{f^{(k)}\}$ *such that* $\theta^{(k)} \geq 0$, $f^{(k)}$ *is bounded below and for all* $k, n \in \mathbb{N}$ *we have* $\theta^{(k+n)} \leq \tau\phi^{(k)}$ *for some fixed* $\tau \in (0,1)$. *Let* $\kappa_f, \delta, \eta \in (0,1)$ *and* $\gamma_s \in (0.5, 1)$ *be given constants. Then for every iteration* $k$, *one of the following holds:*

$$\boldsymbol{\theta - type\ step}: \ \theta^{(k+1)} \leq \tau\phi^{(k)} \text{ and } \phi^{(k+1)} = \left(1-\kappa_f\right)\theta^{(k+1)} + \kappa_f\phi^{(k)}, \tag{A4}$$

$$\boldsymbol{f - type\ step}: \ f^{(k)} - f^{(k+1)} \geq \eta\delta(\theta^{(k)})^{\gamma_s} \text{ and } \phi^{(k+1)} = \phi^{(k)}. \tag{A5}$$

*In both (A4) and (A5) cases, it follows that* $\theta^{(k)} \to 0$ *for* $k \to \infty$.

**Proof.** We analyse two cases, depending on whether there are an infinite number of $\theta - type$ steps:

1. Infinite $\theta - type$ steps. From the update rule (A4), we have:

$$\begin{aligned}
\phi^{(k+1)} &= \left(1-\kappa_f\right)\theta^{(k+1)} + \kappa_f\phi^{(k)} \\
&\leq \left(1-\kappa_f\right)\tau\phi^{(k)} + \kappa_f\phi^{(k)} \\
&\leq (1-(1-\tau)(1-\kappa_f))\phi^{(k)}.
\end{aligned}$$

   Define $q \stackrel{\text{def}}{=} 1-(1-\tau)(1-\kappa_f) \in (0,1)$. Then $\phi^{(k+1)} \leq q\phi^{(k)}$ implies $\phi^{(k)} \to 0$ for $k \to \infty$.

2. Finite $\theta - type$ steps. If only finitely many $\theta - type$ steps occur, there exists an iteration index $k'$ such that for all $k \geq k'$, only $f - type$ steps are taken, satisfying both the switching condition (14) and the Armijo sufficient decrease condition (15). Summing (A5) from $k = k'$ to $k = k' + N$ gives

$$\sum_{k=k'}^{k'+N} (f^{(k)} - f^{(k+1)}) \geq \eta\delta \sum_{k=k'}^{k'+N} (\theta^{(k)})^{\gamma_s},$$

which is

$$f^{(k')} - f^{(k'+N+1)} \geq \eta\delta \sum_{k=k'}^{k'+N} (\theta^{(k)})^{\gamma_s}.$$

Since $f$ is bounded below by assumption AS1, the left-hand side $(f^{(k')} - f^{(k'+N+1)})$ is uniformly bounded above for all $N$. Consequently, the partial sums $\sum_{k=k'}^{k'+N}(\theta^{(k)})^{\gamma_s}$ are bounded for all $N$, which guarantees the convergence of the series $\sum_{k=k'}^{\infty}(\theta^{(k)})^{\gamma_s} < \infty$ (see Chapter 3 in [37] and Theorem 3.2 in [38] for proof). Because $\theta^{(k)} \geq 0$ for all $k$, convergence of the series implies $(\theta^{(k)})^{\gamma_s} \to 0$ as $k \to \infty$. Since $\gamma_s > 0$, we conclude that $\theta^{(k)} \to 0$ as $k \to \infty$.

Thus, in both cases, the sequence of infeasibility measures converges to zero (feasibility), proving that

$$\lim_{k\to\infty} \theta^{(k)} = 0. \tag{A6}$$

### Global Convergence to Stationary Points

We adapt classical NLP convergence arguments for filter-SQP [39][40] and funnel-SQP [33][28] to establish global convergence of the TR funnel algorithm for grey-box optimisation. Detailed lemmas are given next; derivations are omitted where they follow directly from established NLP theory.

In standard NLP algorithms, finite termination occurs when a point $x^{(k)}$ is found such that both the first-order optimality measure and feasibility violation vanish (i.e., $\chi^{(k)} = \theta^{(k)} = 0$). In our grey-box TR funnel method, we additionally have $\sigma^{(k)} \to 0$ ($\sigma^{(k)} \leq \Delta^{(k)}$). Thus, finite termination here means the sequence $\{x^{(k)}\}$ generated by TR funnel grey-box algorithm is finite due to the algorithm stalling within steps 2–4 while reducing the sampling region size and other measures. Under MFCQ, such an iterate $x^{(k)}$ is a first-order critical point. This follows from standard NLP arguments showing that $\chi^{(k)}$, $\theta^{(k)}$, and $\sigma^{(k)}$ all approach zero simultaneously.

**Lemma 1.** (Bound on normal step) *Suppose the TR funnel algorithm is applied on problem (2), and finite termination does not occur. If AS2 and AS6 hold, and $\theta^{(k)} \leq \delta_d$ for some constant $\delta_d > 0$, then there exists a constant $\kappa_{lsc} > 0$, independent of k, such that*

$$\kappa_{lsc}\theta^{(k)} \leq \|d^{(k)}\|. \tag{A7}$$

**Proof.** This follows directly from Lemma 3.1 in [39] for NLP, the same reasoning applies in the grey-box case.

**Lemma 2.** (Convergence to first-order critical point) *Suppose the TR funnel is applied to problem (2), and finite termination does not occur. Let AS1, AS2, AS3, AS4, and AS6 hold. Suppose there exists a subsequence* $\{k_i\} \subset \mathbb{N}$ *such that TRSP is compatible for all* $k_i$*, and*

$$\lim_{i\to\infty} \chi^{(k_i)} = 0, \quad \lim_{i\to\infty} \theta^{(k_i)} = 0, \quad \lim_{i\to\infty} \sigma^{(k_i)} = 0. \qquad (A8)$$

*Then every limit point of the subsequence* $\{x^{(k_i)}\}$ *is a first-order critical point for problem (2).*

**Proof.** Assumption AS2 ensures that the sequence $\{x^{(k)}\}$ is bounded and has at least one limit point. Let $x^*$ be a limit point of $\{x^{(k_i)}\}$. From AS1, all functions involved in problem (2) are twice continuously differentiable. By AS4, MFCQ holds at all feasible limit points.

We define the RM-based criticality measure $\chi'(x, \mathrm{e})$ as:

$$\chi'(x, \mathrm{e}) = |min\, \nabla f(x)^T v| \qquad (A9)$$

$$s.t.\, \nabla h(x)^T v = 0, g(x) + \nabla g(x)^T v \leq 0, v_y - (\nabla r(w) + \mathrm{e})^T v_w = 0, \|v\|_1 \leq 1,$$

where $\mathrm{e}$ is the RM gradient approximation error. If $\chi'(x^*, 0) = \theta(x^*) = 0$, then $x^*$ satisfies the first-order KKT conditions for problem (2).

AS4 is critical in establishing a link between optimality of problem (2) and criticality measure $\chi'$. MFCQ and the twice-continuously differentiable condition (AS1) ensure that the $\chi'$ is a continuous function of limit points [41][42]. Therefore, we say that MFCQ holds at $x^*$ (AS4), and because of Theorem 2.2.1 by Fiacco [43] combined with Theorem 2.3 by Robinson [44], the value $\chi'(x, \mathrm{e})$ is continuous at $(x^*, 0)$.

Let

$$\mathrm{e}^{(k)} := \nabla r^{(k)}\left(x^{(k)} + d^{(k)}\right) - \nabla t\left(x^{(k)} + d^{(k)}\right),$$

so that the observed criticality measure satisfies:

$$\chi^{(k)} = \chi\left(x^{(k)} + d^{(k)}\right) = \chi'\left(x^{(k)} + d^{(k)}, \mathrm{e}^{(k)}\right).$$

From the *κ-fully linear property* AS3 and compatibility of TRSP at $k_i$ (showing that $x^{(k_i)} + d^{(k_i)}$ lies within the trust-region), we know:

$$\left\|\mathrm{e}^{(k_i)}\right\| \leq \kappa_g \sigma^{(k_i)} \to 0, \qquad as\ i \to \infty.$$

From AS6, $\theta^{(k_i)} \to 0$ implies $\left\|d^{(k_i)}\right\| \to 0$. Since $\sigma^{(k_i)} \to 0$, we can now show that

$$\lim_{i\to\infty} \mathrm{e}^{(k_i)} = 0.$$

Finally, combining continuity of $\chi'(x, \mathrm{e})$ at $(x^*, 0)$ with convergence $x^{(k_i)} + d^{(k_i)} \to x^*$ and $\mathrm{e}^{(k_i)} \to 0$, we conclude:

$$\chi'(x^*, 0) = \lim_{i\to\infty} \chi'\left(x^{(k_i)} + d^{(k_i)}, \mathrm{e}^{(k_i)}\right) = \lim_{l\to\infty} \chi^{(k_i)} = 0.$$

Thus, $x^*$ satisfies first-order optimality conditions (as in Definition 2) for problem (2).

Since all $x^{(k_i)}$ are feasible and the feasible region is closed, the limit point $x^*$ is also feasible. By continuity, we have $\theta(x^*) = 0$, confirming first-order stationarity in the original objective.

**Lemma 3.** (Limit of $\theta^{(k)}$ for funnel iterates) *Assume that the TR funnel algorithm is applied to problem (2), and that finite termination does not occur. Let AS1 and AS2 hold, and suppose the set*

$$Z := \{k \geq 0: x^{(k)} \ lies \ in \ the \ funnel\}$$

*is infinite. Then:*

$$\lim_{\substack{k\to\infty \\ k\in Z}} \theta^{(k)} = 0.$$

**Proof.** This follows directly from Theorem 1 (result (A6)). The reduction in sampling region size $\sigma^{(k)}$ to improve RM accuracy (Lemma 2) strengthens the proof of Lemma 3.

**Lemma 4.** (Upper bound on $\theta^{(k)}$) *Assume the TR funnel algorithm is applied to problem (2), and finite termination does not occur. Suppose that AS1, AS2, and AS3 hold, and that the TRSP is compatible at iteration* $k$*. Then, there exists a constant* $\kappa_{ubt} \geq 0$*, independent of* $k$*, such that*

$$\theta^{(k)} \leq \kappa_{ubt}\left(\Delta^{(k)}\right)^{1+\mu}, \tag{A10}$$

$$\theta\left(x^{(k)} + s^{(k)}\right) \leq \kappa_{ubt}\left(\Delta^{(k)}\right)^{2}. \tag{A11}$$

**Proof.** Since TRSP is compatible, using condition (A7) and the trust-region constraint from compatibility problem (8), we say that

$$\kappa_{lsc}\theta^{(k)} \leq \left\|d^{(k)}\right\| \leq \kappa_{\Delta}\Delta^{(k)} min\left[1, \kappa_{\mu}{\Delta^{(k)}}^{\mu}\right], \tag{A12}$$

$$\kappa_{lsc}\theta^{(k)} \leq \left\|d^{(k)}\right\| \leq \kappa_{\Delta}\kappa_{\mu}\left(\Delta^{(k)}\right)^{1+\mu}. \tag{A13}$$

Hence, this proves (A10) for ${\kappa_{ubt}}' = \kappa_{\Delta}\kappa_{\mu}/\kappa_{lsc}$:

$$\theta^{(k)} \leq {\kappa_{ubt}}'\left(\Delta^{(k)}\right)^{1+\mu}. \tag{A14}$$

As TRSP is compatible for iteration $k$, it ensures that $y - r^{(k)}(w) = 0$. From *κ-fully linear* property (3), we have

$$\theta\left(x^{(k)} + s^{(k)}\right) \leq \kappa_h\left(\sigma^{(k)}\right)^2 \leq \kappa_h\left(\Delta^{(k)}\right)^2. \tag{A15}$$

Therefore compatibility (A14) and *κ-fully linear* property (A15) gives us

$$\kappa_{ubt} = max\left[\kappa_h, \kappa_{\Delta}\kappa_{\mu}/\kappa_{lsc}\right],$$

ensuring both (A10) and (A11) hold.

**Lemma 5.** (Lower bound on objective improvement) *Assume that the TR funnel algorithm is applied to problem (2) and AS1, AS2, AS5, and AS6 hold. Suppose finite termination does not occur, that* $k \notin \mathcal{R}$ *(*$\mathcal{R}$ *is a set of restoration iterations), and that* $\chi^{(k)} \geq \epsilon > 0$*. Furthermore,*

$$\Delta^{(k)} \leq min\left[\left(\frac{\epsilon}{\kappa_{umh}}\right), \left(\frac{2\kappa_{ubg}}{\kappa_{umh}\kappa_{\Delta}\kappa_{\mu}}\right)^{\frac{1}{1+\mu}}, \left(\frac{\kappa_{tmd}\epsilon}{4\kappa_{ubg}\kappa_{\Delta}\kappa_{\mu}}\right)^{\frac{1}{\mu}}\right] \stackrel{\text{def}}{=} \delta_m \tag{A16}$$

*where* $\kappa_{ubg} = \max_{x\in X}\|\nabla f(x)\|$ *and* $\kappa_{umh} = max\left[\max_{k}\{\beta^{(k)}\}, \max_{x\in X}\|\nabla^2 f(x)\|\right]$. *Then,*

$$f\left(x^{(k)}\right) - f\left(x^{(k)} + s^{(k)}\right) \geq \frac{1}{2}\kappa_{tmd}\epsilon\Delta^{(k)}. \tag{A17}$$

**Proof.** We prove Lemma 5 using the classical Cauchy decrease lower bound from trust-region literature [38]. By AS5, trial step $s^{(k)}$ satisfies fraction of Cauchy decrease on the actual objective:

$$f\left(x^{(k)}\right) - f\left(x^{(k)} + s^{(k)}\right) \geq \kappa_{tmd}\left\|\nabla f\left(x^{(k)}\right)\right\| min\left[\frac{\left\|\nabla f\left(x^{(k)}\right)\right\|}{\left\|\nabla^2 f\left(x^{(k)}\right)\right\|}, \Delta^{(k)}\right]. \tag{A18}$$

Since $\left\|\nabla f\left(x^{(k)}\right)\right\| \geq \chi^{(k)} \geq \epsilon$ and $\|\nabla^2 f(x)\| \leq \kappa_{umh}$, then

$$\frac{\left\|\nabla f\left(x^{(k)}\right)\right\|}{\left\|\nabla^2 f\left(x^{(k)}\right)\right\|} \geq \frac{\epsilon}{\kappa_{umh}}.$$

Therefore, if $\Delta^{(k)} \leq \epsilon/\kappa_{umh}$, the min operator in (A18) on the right-hand side of the inequality becomes $\Delta^{(k)}$, and we conclude:

$$f\left(x^{(k)}\right) - f\left(x^{(k)} + s^{(k)}\right) \geq \kappa_{tmd}\left\|\nabla f\left(x^{(k)}\right)\right\|\Delta^{(k)} \geq \kappa_{tmd}\epsilon\Delta^{(k)}. \tag{A19}$$

Inequality (A16) imposes slightly stronger (conservative) bounds so that, after tightening constants, we obtain (A17). We keep $\Delta^{(k)}$ in Lemmas 5-8 because the Cauchy-type lower bound concerns globalisation of the algorithm rather than RM accuracy.

**Lemma 6.** ($f - type$ for small trust-region radius) *Assume that the TR funnel algorithm is applied to problem (2) and AS1, AS2, AS3, AS5, and AS6 hold. Suppose finite termination does not occur, that* $k \notin \mathcal{R}$ *(iterations don't belong to the restoration phase), (A7) holds, and that* $\chi^{(k)} \geq \epsilon > 0$. *Furthermore,*

$$\Delta^{(k)} \leq min\left[\mathrm{e}_m, \left(\frac{\kappa_{tmd}\epsilon}{2\kappa_{\theta}{\kappa_{ubt}}^{\gamma_s}}\right)^{\frac{1}{\gamma_s(1+\mu)-1}}\right] \stackrel{\text{def}}{=} \delta_f. \tag{A20}$$

Then,

$$f^{(k)} - f^{(k+1)} \geq \eta\delta(\theta^{(k)})^{\gamma_s}.$$

**Proof.** Using mathematical operations (e.g., exponent simplification and re-arranging) on (A10), Lemma 4, (A20), and Lemma 5, we get

$$\eta\delta(\theta^{(k)})^{\gamma_s} \leq \eta\delta\left(\kappa_{ubt}\left(\Delta^{(k)}\right)^{1+\mu}\right)^{\gamma_s} \leq \frac{1}{2}\kappa_{tmd}\epsilon\Delta^{(k)} \leq f\left(x^{(k)}\right) - f\left(x^{(k)} + s^{(k)}\right).$$

**Lemma 7.** (Small infeasibility measure or trust-region radius guarantee funnel acceptance) *Assume that the TR funnel algorithm is applied to problem (2) and AS1, AS2, AS3, AS5, and AS6 hold. Suppose finite termination does not occur,* $k \notin \mathcal{R}$,*and* $\chi^{(k)} \geq \epsilon > 0$, *and*

$$\theta^{(k)} \leq \kappa_{ubt}\left(\Delta^{(k)}\right)^{1+\mu} \qquad (:= \delta_\theta), \tag{A21}$$

*Then*

$$f\left(x^{(k)} + s^{(k)}\right) \leq f\left(x^{(k)}\right) - \eta\delta(\theta^{(k)})^{\gamma_s}$$

**Proof.** From (A10), we have

$$\Delta^{(k)} \geq \left(\frac{\theta^{(k)}}{\kappa_{ubt}}\right)^{\frac{1}{1+\mu}}. \tag{A22}$$

Using (A21) and (A22), implies that $\Delta^{(k)}$ satisfies $\Delta^{(k)} \leq \delta_m$ (A20). Since all the assumptions of Lemma 6 hold, it directly proves Lemma 7.

**Lemma 8.** (No restoration iteration when infeasibility measure and trust-region radius are small) *Assume that the TR funnel algorithm is applied to problem (2) and AS1, AS2, AS3, AS5, and AS6 hold. Suppose finite termination does not occur and for some iteration* $k > 0$, $\chi^{(k)} \geq \epsilon > 0$, *and*

$$\Delta^{(k)} \leq min\left[\left(\frac{1}{\kappa_\mu}\right)^{\frac{1}{\mu}}, \left(\frac{\gamma_c^2 c \kappa_\Delta \kappa_\mu}{\kappa_{usc}\kappa_{ubt}}\right)^{\frac{1}{(1-\mu)}}\right] \stackrel{\text{def}}{=} \delta_R, \tag{A23}$$

*where* $c \in (0, 1)$. *Further, suppose*

$$\theta^{(k)} \leq min[\delta_\theta, \delta_d]. \tag{A24}$$

*Then* $k \notin \mathcal{R}$.

**Proof.** By the stated assumptions and the bound on $\theta^{(k)} \leq \delta_d$ (A24), Lemma 1 yields

$$\kappa_{lsc}\theta^{(k)} \leq \left\|d^{(k)}\right\| \leq \kappa_{usc}\theta^{(k)}. \tag{A25}$$

From $\theta^{(k)} \leq \delta_\theta = \kappa_{ubt}\left(\Delta^{(k)}\right)^{1+\mu}$ and Lemma 7, we also have the funnel acceptance condition:

$$f\left(x^{(k)} + s^{(k)}\right) \leq f\left(x^{(k)}\right) - \eta\delta(\theta^{(k)})^{\gamma_s}$$

Suppose, for contradiction, that iteration $k \in \mathcal{R}$. Then, by the definition of a restoration iteration in the algorithm, we must have

$$\left\|d^{(k)}\right\| > \kappa_\Delta \kappa_\mu (\Delta^{(k)})^{1+\mu}. \tag{A26}$$

Since $\Delta^{(k)} \leq \left(\frac{1}{\kappa_\mu}\right)^{\frac{1}{\mu}}$ (from $\delta_R$), the right-hand side of (A26) is positive.

In first case, consider $k-1$ iteration unsuccessful. From the unsuccessful step update, we get $\Delta^{(k-1)} \geq \gamma_c \Delta^{(k)}$. Lemma 4 then gives:

$$\theta\left(x^{(k-1)} + s^{(k-1)}\right) \leq \kappa_{ubt}\left(\Delta^{(k-1)}\right)^2 \leq \frac{\kappa_{ubt}}{\gamma_c^2}\left(\Delta^{(k)}\right)^2. \tag{A27}$$

Using (A25), (A26) and (A27), we get

$$\kappa_\Delta \kappa_\mu (\Delta^{(k)})^{1+\mu} < \left\|d^{(k)}\right\| \leq \kappa_{usc}\theta^{(k)} \leq \kappa_{usc}\kappa_{ubt}\left(\Delta^{(k-1)}\right)^2 \leq \frac{\kappa_{usc}\kappa_{ubt}}{\gamma_c^2}\left(\Delta^{(k)}\right)^2,$$

or

$$\kappa_\Delta \kappa_\mu (\Delta^{(k)})^{1+\mu} < \left\|d^{(k)}\right\| \leq \kappa_{usc}\theta^{(k)} \leq \frac{\kappa_{usc}\kappa_{ubt}}{\gamma_c^2}\left(\Delta^{(k)}\right)^2. \tag{A28}$$

Dividing (A28) by $(\Delta^{(k)})^{1+\mu} > 0$ gives:

$$\kappa_\Delta \kappa_\mu < \frac{\kappa_{usc}\kappa_{ubt}}{\gamma_c^2}\left(\Delta^{(k)}\right)^{1-\mu}. \tag{A29}$$

Since $0 < \mu < 1$, inequality (A29) implies:

$$(\Delta^{(k)})^{1-\mu} > \frac{\gamma_c^2 \kappa_\Delta \kappa_\mu}{\kappa_{usc}\kappa_{ubt}},$$

which contradicts $\Delta^{(k)} \leq \delta_R$.

In this second case, we consider $k-1$ iteration is successful. A similar argument applies, using the fact that in a successful step, $\Delta^{(k)} \geq \gamma_e \Delta^{(k-1)}$ with $\gamma_e > 1$ and (A21), leading again to a contradiction with $\Delta^{(k)} \leq \delta_R$.

Therefore, $k \in \mathcal{R}$ is impossible under the stated conditions, completing the proof.

**Lemma 9.** (Convergence along infinite funnel subsequence) *Assume the TR funnel algorithm is applied to problem (2) and that finite termination does not occur. Let AS1, AS2, AS3, AS5, and AS6 hold. Let*

$$\mathcal{Z} = \{k \geq 0 \mid k \notin R\}$$

*be the (infinite) set of iterations at which the restoration is not invoked. Then there exists a subsequence* $\{k_l\} \subseteq \mathcal{Z}$ *such that*

$$\lim_{l\to\infty} \theta^{(k_l)} = 0, \tag{A30}$$

$$\lim_{l\to\infty} \chi^{(k_l)} = 0, \tag{A31}$$

$$\lim_{k\to\infty} \sigma^{(k_l)} = 0. \tag{A32}$$

**Proof.** The first limit (A30) follows directly from Lemma 3 and convergence to feasibility proof: since $\mathcal{Z}$ is infinite, we can select a subsequence $\{k_l\} \subseteq \mathcal{Z}$ for which $\theta^{(k_l)} \to 0$.

Assume, for contradiction, that $\chi^{(k_l)} \nrightarrow 0$. Then, there exists $\epsilon_2 > 0$ and a further subsequence $\{k_m\} \subseteq \{k_l\}$ with $\chi^{(k_m)} \geq \epsilon_2, \forall m$. Consider the trust-region radii along $\{k_m\}$. Two cases are possible:

- Case 1: There exists $\iota > 0$ and an infinite sub-subsequence with $\Delta^{(k_m)} \geq \iota$. For each accepted (successful) iteration in this sub-subsequence, Lemma 5 yields

$$f\left(x^{(k_m)}\right) - f\left(x^{(k_m)} + s^{(k_m)}\right) \geq \frac{1}{2}\kappa_{tmd}\epsilon_2\Delta^{(k_m)} \geq c_2 > 0,$$

with $c_2 := \frac{1}{2}\kappa_{tmd}\epsilon_2\iota$. Since $f$ is bounded below and infinitely many such accepted decreases cannot occur, Case 1 is impossible.

- Case 2: Assume $\Delta^{(k_m)} \rightarrow 0$. Because $\theta^{(k_m)} \rightarrow 0$ as $m \rightarrow \infty$, for sufficiently large $m$, $\Delta^{(k_m)} < \delta_f$. Lemma 6 then guarantees that, since $\chi^{(k_m)} \geq \epsilon_2$, $f - type$ (accepted objective-improving) steps occur for all large $m$. But infinitely many $f - type$ successes would prevent $\Delta^{(k_m)} \rightarrow 0$ (and $\sigma^{(k_m)} \rightarrow 0$, as $\sigma$ remains unchanged), contradicting the assumption of Case 2.

Both cases lead to contradictions; hence $\chi^{(k_l)} \rightarrow 0$ as $k_l \rightarrow \infty$, proving (A31).

Finally, if $\sigma^{(k_l)} \nrightarrow 0$, then there exists $u > 0$ and an infinite sub-subsequence $\{k_s\} \subseteq \{k_l\}$ with $\sigma^{(k_s)} \geq u$ for all $s$. If infinitely many of these indices correspond to successful iterations in terms of objective improvement, Lemma 5 yields a uniform positive reduction of $f$ for each, contradicting boundedness of $f$. Hence, only finitely many can be accepted; therefore both $\Delta^{(k)}$ and $\sigma^{(k)}$ decrease along an infinite subsequence, enabling extraction of a sub-subsequence with $\sigma^{(k_l)} \rightarrow 0$ as $k_l \rightarrow \infty$, proving (A32).

Combining the subsequence extraction (A32) with (A30) and (A31) yields the desired subsequence $\{k_l\} \subseteq \mathcal{Z}$ with $\theta^{(k_l)} \rightarrow 0$, $\chi^{(k_l)} \rightarrow 0$ and $\sigma^{(k_l)} \rightarrow 0$ as $k_l \rightarrow \infty$.

**Lemma 10.** ($\theta \rightarrow 0$ for finite funnel subsequence) *Assume the TR funnel algorithm is applied to problem (2) and that finite termination does not occur. Let AS1, AS2, AS3, AS5, and AS6 hold, and suppose the set*

$$Z := \left\{k \geq 0 : x^{(k)} \text{ lies in the funnel}\right\}$$

*is finite. Then,*

$$\lim_{k \rightarrow \infty} \theta^{(k)} = 0. \qquad \textit{(A33)}$$

**Proof.** Let $k_0 > 0$ as the last index with $x^{(k_0-1)}$ that updates the funnel width ($\theta - type$). Then for all $k \geq k_0$, every successful iteration is an $f - type$ step.

If the set of successful iterates with $k \geq k_0$ were finite, then $\Delta^{(k)} \rightarrow 0$, and Lemma 7 would imply finite termination (contradicting our hypothesis). Hence, there are infinitely many $f - type$ successful iterations after $k_0$.

For each such successful iteration $k \geq k_0$ we have,

$$f\left(x^{(k)}\right) - f\left(x^{(k+1)}\right) \geq \eta\delta(\theta(x^{(k)}))^{\gamma_s} \geq 0. \qquad \text{(A34)}$$

By AS1 and AS2, $f$ is bounded below and nonincreasing for $k \geq k_0$, so

$$\lim_{k\to\infty} f\left(x^{(k)}\right) - f\left(x^{(k+1)}\right) = 0. \tag{A35}$$

Combining (A35) with (A34), we conclude that $\theta^{(k)} \to 0$ as $k \to \infty$, proving (A33).

**Lemma 11.** (for a finite funnel sequence, $\Delta^{(k)}$ is bounded away from zero if $\chi^{(k)}$ bounded away from zero) *Assume the TR funnel algorithm is applied to problem (2) and that finite termination does not occur. Let AS1, AS2, AS3, AS5, and AS6 hold, and suppose the set*

$$Z := \left\{k \geq 0: x^{(k)}\ lies\ in\ the\ funnel\right\}$$

*is finite. If* $\chi^{(k)} \geq \epsilon > 0$ *for all* $k \geq k_0$ *(where* $k_0$ *is defined as in the lemma 10), then there exists a* $\Delta_{min}> 0$ *such that*

$$\Delta^{(k)} \geq \Delta_{min}$$

*for all* $k$.

**Proof.** Since $Z$ is finite, all successful iterations $k \geq k_0$ yield $f - type$ steps, and no restoration iteration occurs. Assumptions AS1–AS3, AS5, and AS6 guarantee that $f$ is sufficiently smooth, the criticality measure $\chi^{(k)}$ is well defined, and the trust-region updates satisfy the standard conditions of Lemma 3.12 in [39]. The hypothesis $\chi^{(k)} \geq \epsilon > 0$ for all $k \geq k_0$ ensures that the lower-bound condition on the criticality measure in Lemma 3.12 [39] is satisfied.

Applying Lemma 3.12 in [39] therefore yields the existence of $\Delta_{min}> 0$ such that $\Delta^{(k)} \geq \Delta_{min}\ \forall k$.

**Lemma 12.** (Convergence of $\chi$) *Assume the TR funnel algorithm is applied to problem (2) and that finite termination does not occur. Let AS1, AS2, AS3, AS5, and AS6 hold, and suppose the set*

$$Z := \left\{k \geq 0: x^{(k)}\ lies\ in\ the\ funnel\right\}$$

*is finite. Then,*

$$\lim_{k\to\infty} \inf \chi^{(k)} = 0. \tag{A36}$$

**Proof.** From Lemma 10, we have $\theta^{(k)} \to 0$, and by AS6 this implies $\left\|d^{(k)}\right\| \to 0$. The convergence of $f$ along successful iterates, as in (A35), yields

$$\lim_{\substack{k\to\infty \\ k\in S}} \left[f\left(x^{(k)} + d^{(k)}\right) - f\left(x^{(k)} + s^{(k)}\right)\right] = 0. \tag{A37}$$

Suppose, for contradiction, that $\chi^{(k)} \geq \epsilon > 0$ for all $k \geq k_0$. By Lemma 11 and applying the fraction of Cauchy decrease condition then gives, for all $k \geq k_0$,

$$f\left(x^{(k)} + d^{(k)}\right) - f\left(x^{(k)} + s^{(k)}\right) \geq \kappa_{tmd} \epsilon min \left[\frac{\epsilon}{\kappa_{umh}}, \Delta_{min}\right] \geq c_* > 0. \tag{A38}$$

This (A38) contradicts (A37), which states that the left-hand side tends to zero along successful iterates.

Therefore, $\chi^{(k)}$ cannot remain bounded away from zero and (A36) follows.

Finally, the convergence proof is summarised in the following theorem.

**Theorem 2.** *Assume the TR funnel algorithm is applied to problem (2) and that finite termination does not occur. Let AS1, AS2, AS3, AS4, AS5, and AS6 hold. Then either the restoration phase terminates unsuccessfully at a stationary point after exhausting the computational budget or there exists a subsequence* $\{k_l\}$ *for which*

$$\lim_{l\to\infty} x^{(k_l)} = x^* \tag{A39}$$

*where* $x^*$ *is a first-order KKT point of problem (2).*

**Proof.** Suppose that the restoration phase always terminates successfully. From Lemmas 9, 10, and 12, there exists a subsequence $\{k_l\}$ such that

$$\lim_{l\to\infty} \chi^{(k_l)} = 0, \quad \lim_{l\to\infty} \theta^{(k_l)} = 0.$$

By step 3(c) in the TR funnel algorithm and Lemma 2, we have

$$\lim_{l\to\infty} \sigma^{(k_l)} = 0.$$

Finally, applying Lemma 2 yields the conclusion that $x^*$ is a first-order KKT point of problem (2).

# Annexe A: Mathematical Formulations of Optimisation Problems

## Colville Problem (AN1)

We consider a modified version of the Colville function reformulated as a constrained grey-box optimisation problem (A1) with variable dimensions $(n_w, n_y, n_z) = (4,4,1)$:

$$\min_{w,y,z} \quad 5.3578w_3^2 + y_1, \tag{AN1}$$

$$\begin{aligned}
s.t. \quad & g_1(w,y,z) := y_2 - 0.0000734w_1w_4 - 1 \leq 0, \\
& g_2(w,y,z) := 0.000853007w_2z_5 + 0.00009395w_1w_4 - 0.00033085w_3z_5 - 1 \leq 0, \\
& g_3(w,y,z) := y_4 - 0.30586(w_2z_5)^{-1}w_3^2 - 1 \leq 0, \\
& g_4(w,y,z) := 0.00024186w_2z_5 + 0.00010159w_1w_2 + 0.00007379w_3^2 - 1 \leq 0, \\
& g_5(w,y,z) := y_3 - 0.40584(z_5)^{-1}w_4 - 1 \leq 0, \\
& g_6(w,y,z) := 0.00029955w_3z_5 + 0.00007992w_1w_3 + 0.00012157w_3w_4 - 1 \leq 0, \\
& y_1 = d_1(w) := 0.8357w_1w_4 + 37.2392w_1, \\
& y_2 = d_2(w) := 0.00002584w_3w_4 + 0.00006663w_2w_4, \\
& y_3 = d_3(w) := 2275.1327(w_3w_4)^{-1} + 0.2668(w_4)^{-1}w_1, \\
& y_4 = d_4(w) := 1330.3294(w_2w_4)^{-1} + 0.42(w_4)^{-1}w_1, \\
& 78 \leq w_1 \leq 102, \\
& 33 \leq w_2 \leq 45, \\
& 27 \leq w_i \leq 45, \forall i \in \{3,4\}, \\
& 27 \leq z_5 \leq 45.
\end{aligned}$$

## Himmelblau Problem (AN2)

We propose a modified formulation of the Himmelblau problem that can be formulated as a constrained grey-box optimisation problem (A2) with number of variables $(n_w, n_y, n_z) = (3,2,5)$, as follows

$$\min_{w,y,z} \quad 5.3578547y_1 + 0.8356891z_1w_5 + 37.2932239z_1 - 40792.141, \tag{AN2}$$

$$\begin{aligned}
s.t. \quad & h_1(w,y,z) := z_6 = 85.334407 + 0.0056858y_2 + 0.00026z_1z_4 - 0.0022053w_3w_5, \\
& h_2(w,y,z) := z_7 = 80.51249 + 0.0071317y_2 + 0.0029955z_1w_2 - 0.0021813{w_3}^2, \\
& h_3(w,y,z) := z_8 = 9.300961 + 0.0047026w_3w_5 + 0.0012547z_1w_3 - 0.0019085w_3z_4, \\
& y_1 = d_1(w) := w_3^2, \\
& y_2 = d_2(w) := w_2w_5, \\
& 78 \leq z_1 \leq 102,
\end{aligned}$$

$33 \leq w_2 \leq 45,$

$27 \leq w_i \leq 45, \forall i \in \{3,5\},$

$27 \leq z_4 \leq 45,$

$0 \leq z_6 \leq 92,$

$90 \leq z_7 \leq 110,$

$20 \leq z_8 \leq 25.$

### Loeppky Problem (AN3)

We propose a modified formulation of the Loeppky function that can be formulated as a constrained grey-box optimisation problem (A3) with number of variables $(n_w, n_y, n_z) = (3,1,4)$, as follows

$$\min_{w,y,z} \quad 6w_1 + 4w_2 + 5.5w_3 + y_1 + 1.4w_2w_3 + z_4 + 0.5z_5 + 0.2z_6 + 0.1z_7, \tag{AN3}$$

$$s.t. \quad y_1 = d_1(w) := 3w_1w_2 + 2.2w_1w_3,$$

$$0 \leq w_i, z_j \leq 1, \forall i \in \{1,2,3\}, \forall j \in \{4,5,6,7\}.$$

### Wing Weight Problem (AN4)

We propose a modified formulation of the wing weight function (modelling a light aircraft wing) that can be formulated as a constrained grey-box optimisation problem (A4) with number of variables $(n_w, n_y, n_z) = (2,1,8)$, as follows

$$\min_{w,y,z} \quad 0.036{w_1}^{0.758}{z_2}^{0.0035}\left(\frac{z_3}{\cos^2(z_4)}\right)^{0.6}{z_5}^{0.006}{z_6}^{0.04}\left(\frac{100z_7}{\cos(z_4)}\right)^{-0.3}(z_8z_9)^{0.49} + y_1, \tag{AN4}$$

$$s.t. \quad y_1 = d_1(w) := w_1w_2,$$

$$150 \leq w_1 \leq 200,$$

$$220 \leq z_2 \leq 300,$$

$$6 \leq z_3 \leq 10,$$

$$-10 \leq z_4 \leq 10,$$

$$16 \leq z_5 \leq 45,$$

$$0.5 \leq z_6 \leq 1,$$

$$0.08 \leq z_7 \leq 0.18,$$

$$2.5 \leq z_8 \leq 6,$$

$$1700 \leq z_9 \leq 2500,$$

$$0.025 \leq w_{10} \leq 0.08.$$

### Welded Beam Problem (AN5)

A modified formulation to design a welded beam with the least cost can be presented as a constrained grey-box optimisation problem (A5) with number of variables $(n_w, n_y, n_z) = (4,1,0)$, as follows

$$\min_{w,y,z} \quad y_1, \tag{AN5}$$

$$
\begin{aligned}
s.t. \quad & g_1(w,y,z) := \tau(w,y,z) - 13000 \le 0,\\
& g_2(w,y,z) := \sigma(w,y,z) - 30000 \le 0,\\
& g_3(w,y,z) := w_1 - w_4 \le 0,\\
& g_4(w,y,z) := 0.10476{w_1}^2 + 0.04811 w_3 w_4 (14 + w_2) - 5 \le 0,\\
& g_5(w,y,z) := 0.125 - w_1 \le 0,\\
& g_6(w,y,z) := \delta(w,y,z) - 0.25 \le 0,\\
& g_7(w,y,z) := 60000 - P_c(w,y,z) \le 0,\\
& y_1 = d_1(w) := 1.10471 {w_1}^2 w_2 + 0.04811 w_3 w_4 (14 - w_2),
\end{aligned}
$$

where:

$$\tau(w,y,z) = \sqrt{(\tau')^2 + 2\tau'\tau''\frac{w_2}{2R} + (\tau'')^2},$$

$$\tau' = \frac{6000}{\sqrt{2} w_1 w_2},$$

$$\tau'' = \frac{MR}{J},$$

$$M = 6000\left(14 + \frac{w_2}{2}\right),$$

$$R = \sqrt{\frac{{w_2}^2}{4} + \left(\frac{w_1 + w_3}{2}\right)^2},$$

$$J = 2\left\{\sqrt{2}\, w_1 w_2 \left[\frac{{w_2}^2}{12} + \left(\frac{w_1 + w_3}{2}\right)^2\right]\right\},$$

$$\sigma(w,y,z) = \frac{504000}{w_4 {w_3}^2},$$

$$\delta(w,y,z) = \frac{2.1952}{{w_3}^3 w_4},$$

$$P_c(w,y,z) = 64746.022(1 - 0.0282346 w_3) w_3 {w_4}^3,$$

$$0.125 \le w_1 \le 5,$$

$$0.1 \le w_2 \le 10,$$

$$0.1 \le w_3 \le 10,$$

$$0.1 \le w_4 \le 5.$$

### Williams-Otto Process (AN6)

A modified formulation of the Williams-Otto process can be presented as a constrained grey-box optimisation problem (A6), where the reactor is considered a black-box, with number of variables $(n_w, n_y, n_z) = (6,3,21)$, as follows

$$\min_{w,y,z} \quad ROI: 100(2207F_P + 50F_{purge} - 168F_A - 252F_B - 2.22F_{eff}^{sum} - 84F_G - 60V\rho)/(600V\rho), \tag{AN6}$$

$$
\begin{aligned}
s.t. \quad & h_1(w,y,z) := F_{eff}^A = F_A + F_R^A - y_1,\\
& h_2(w,y,z) := F_{eff}^B = F_B + F_R^B - (y_1 + y_2),\\
& h_3(w,y,z) := F_{eff}^C = F_R^C + 2r_1 - 2y_2 - y_3,\\
& h_4(w,y,z) := F_{eff}^E = F_R^E + 2y_2,\\
& h_5(w,y,z) := F_{eff}^P = 0.1F_R^E + y_2 - 0.5y_3,\\
& h_6(w,y,z) := F_{eff}^G = 1.5y_3,\\
& h_7(w,y,z) := F_{eff}^{sum} = \sum_j F_{eff}^j, \;\; j \in \{A,B,C,E,P,G\},\\
& h_8(w,y,z) := F_G = \; F_{eff}^G,\\
& h_9(w,y,z) := F_{purge} = \; \eta(F_{eff}^A + F_{eff}^B + F_{eff}^C + 1.1F_{eff}^E),\\
& h_{10}(w,y,z) := F_R^j = (1-\eta)F_{eff}^j, \qquad j \in \{A,B,C,E,P,G\},\\
& y_1 = d_1(w) := k_1 x_A x_B V\rho = a_1 e^{-120/T} x_A x_B V\rho,\\
& y_2 = d_2(w) := k_2 x_B x_C V\rho = a_2 e^{-150/T} x_B x_C V\rho,\\
& y_3 = d_3(w) := k_3 x_P x_C V\rho = a_3 e^{-200/T} x_P x_C V\rho,
\end{aligned}
$$

where

$V \in [0.03, 0.1)$,

$T \in [5.8, 6.8]$,

$F_P \in [0, 4.763]$,

$F_{purge}\, , F_G, F_{eff}^j \geq 0$,

$F_A, F_B \geq 1$,

$a_1 = 5.9755 \times 10^9$,

$a_2 = 2.5962 \times 10^{12}$,

$a_3 = 9.6283 \times 10^{15}$,

$\rho = 50$.

### Biomass-based Hydrogen Production (AN7)

We consider a process systems optimisation problem for minimising the capital (CAPEX) and operating (OPEX) expenditures of a biomass gasification hydrogen production plant. This techno-economic grey-box optimisation problem has $(n_w, n_y, n_z) = (8,4,71)$ variables and 69 constraints. The complete code of the mathematical formulation of this problem is available in the GitHub repository [35].